\documentclass[12pt,leqno]{article}
\usepackage{amssymb}
\usepackage{amsfonts}
\usepackage{graphicx}
\newcommand{\truc}{\;|\!\!\!\!\Rightarrow}
\begin{document}

\begin{center}
{\Huge
cinq

conf\'erences

sur

l'ind\'ecidabilit\'e

}
{\large \vspace{4cm}

organis\'ees en 1982 \`a l'Ecole Nationale des Ponts et Chauss\'ees

dans le cadre du programme

Sciences, Techniques et Soci\'et\'e}

\vspace{5cm}

{\Large Nicolas Bouleau,  JeanYves Girard,  Alain Louveau}
 
\end{center}
\newpage
\begin{center}\Large
{\Large\bf Sommaire}
\end{center}

\vspace{2cm}

	\hfill Pages\\

{\bf Introduction}	\hfill 3\\

{\bf La formalisation des math\'ematiques},	\hfill 6

Nicolas Bouleau\\

{\bf Sur le th\'eor\`eme d'incompl\'etude de G\"{o}del}

     {\bf - jusqu'\`a G\"{o}del}	\hfill 15

             {\bf  - apr\`es G\"{o}del},	\hfill 24

Jean-Yves Girard\\

{\bf Ind\'ecidabilit\'e de l'hypoth\`ese du continu},	\hfill 		     35

Alain Louveau\\

{\bf Sur la calculabilit\'e effective, exemples},	\hfill  43

Nicolas Bouleau\\

{\bf Bibliographie}	\hfill 55
 
\newpage

\vspace{2cm}

\noindent{\LARGE\bf Introduction}\\

{\tt
Entre les scientifiques purs, chercheurs et universitaires, et l'habitant de nos cit\'es dans sa vie quotidienne, existe une cat\'egorie de personnes qui utilisent la science sous diverses formes et la relient \`a  la vie sociale. Appelons ceux-ci les ing\'enieurs, ce sont eux qui font que les trains circulent, que les ponts tiennent, que les avions sont guid\'es. Un tel ing\'enieur poss\`ede un savoir efficace, une th\'eoÂrie pertinente qui permet de d\'eboucher effectivement sur les actions \`a entreprendre : conception, dimensionnement, v\'erification. Que cette th\'eorie ne soit qu'une approximation, l'ing\'enieur le sait, il n'a pas pris en compte la structure fine de la mati\`ere et ses contr\^oles ne sont que statistiques, mais il ma\^{\i}trise les ordres de grandeurs et n\'eglige ce qui est n\'egligeable.\\

Y aurait-il des questions ind\'ecidables ? Notre ing\'enieur aura tendance \`a consid\'erer que ce sont l\`a des \'elucubrations purement sp\'eculatives sans incidences sur la vie r\'eelle qu'il conna\^{\i}t et qui est bien d\'ecrite par les m\'ethodes qu'il utilise.\\

Il n'a raison qu'\`a demi. D\`es qu'on pense aux enjeux, \`a l'int\'er\^et g\'en\'eral et \`a celui de chacun, plusieurs lectures de la r\'ealit\'e apparaissent. Ces conf\'erences sur l'ind\'ecidabilit\'e, peuvent le concerner par analogie entre sa situation et celle des math\'ematiciens.\\

En effet, le plus souvent, le math\'ematicien reste install\'e dans son cadre conceptuel habituel qui est aujourdÕhui celui de la th\'eorie des ensembles et ne voit pas l'int\'er\^et d'en sortir. Jean Dieudonn\'e \'ecrit \`a propos des logiciens "Nous, les math\'ematiciens, comment voyons-nous leur travail ? Eh bien, d'une part, ils explorent les possibilit\'es de notre syst\`eme logique celui avec lequel nous travaillons, Zermelo-Fraenkel, d'un autre c\^ot\'e - et cela nous int\'eresse beaucoup moins - ils \'elaborent et explorent des quantit\'es d'autres syst\`emes logiques. Alors, quand on vient nous parler de la logique du premier et du deuxi\`eme ordre, de fonctions r\'ecursives et de mod\`eles, th\'eories tr\`es gentilles et tr\`es belles qui ont obtenu des r\'esultats remarquables, nous, math\'ematiciens, nous ne voyons aucune objection \`a ce qu'on s'en occupe, mais cela nous laisse enti\`erement froids"\footnote{ {\it Penser les mathŽmatiques}, Seuil 1982.}. Cette attitude sereine semble \^etre une position pertinente. Cependant, lorsqu'un math\'ematicien a rencontr\'e dans la poursuite naturelle de ses recherches un probl\`eme difficile, et que ce probl\`eme apr\`es des tentatives infructueuses de d\'eÂmonstration et de r\'efutation s'av\`ere ind\'ecidable, alors il est convaincu qu'on ne peut laisser la logique de c\^ot\'e avec ses fonctions r\'ecursives et ses mod\`eles, il change de camp, s'int\'eÂresse \`a des syst\`emes syntaxiques ou s\'emantiques diff\'erents (comme celui de l'analyse non standard qui trouve maintenant des applications en m\'ecanique, en physique et en calcul stochastique) comme \`a autant d'outils possibles.\\

De m\^eme l'ing\'enieur ne peut me semble-t-il aujourdÕhui, ne pas se pr\'eoccuper du contexte de son activit\'e, des approches diff\'erentes, du caract\`ere relatif de sa mod\'elisation favorite.\\

Quelques mots, en outre, doivent \^etre dits sur l'int\'er\^et intrins\`eque du sujet. Ce que dit A.Shenitzer \`a propos des \'etudiants en math\'ematiques "Speaking of the deductive method, it is a sad reflexion on the intellectual level of mathematical education that, unless he takes courses in logic, the mathematics student may get his degree without having heard about G\"{o}del or about his monumental discovery of the intrinsic limitations of deductive method, a discovery widely regarded as one of the greatest intellectual accomplishements of the 20th century"\footnote{``Teaching mathematics" in {\it Mathematics tomorrow} L.A. Steen , editor Springer 1981.} s'applique aussi bien \`a l'\'el\`eve-ing\'e\-nieur. C'est tout le probl\`eme de la vulgarisation qui se pose alors. La logique math\'ematique est un domaine tr\`es \'esot\'erique. Que font les logiciens ? Je peux apporter un modeste t\'emoignage, je suis all\'e me promener par l\`a-bas et je puis dire qu'ils font des choses difficiles mais passionnantes, certai\-nement pas anodines, des questions centrales sont abord\'ees, les m\'ethodes et les concepts sont forts et \'eclairants, \`a tel point qu'une certaine philosophie de la repr\'esentation para\^{\i}t bien vaine \`a c\^ot\'e. C'est d'ailleurs ce qui explique que les logiciens r\'epugnent \`a vulgariser leurs travaux. M\^eme si leur d\'emarche est tr\`es abstraite, ils ont le sentiment de tenir en main du r\'eel bien plus que certains verbiages philosophiques et que les enjeux des difficult\'es qu'ils surmontent appara\^{\i}tront un jour ou l'autre. \\

La question de la vulgarisation se pose avec d'autant plus de force. Je pense que les lecteurs trouveront de l'int\'er\^et \`a lire dans les expos\'es de Jean-Yves GIRARD et Alain LOUVEAU une sorte de philosophie vivante des math\'ematiques qui est simplement la vision fra\^{\i}che de chercheurs de haut niveau dans leur propre champ.\\

\hfill Nicolas Bouleau
\newpage
\noindent{\LARGE\bf La formalisation des math\'ematiques}

\noindent Nicolas Bouleau\\

Cette s\'eance a pour but de vous pr\'esenter quelques notions introductives vous permettant de suivre plus ais\'ement les prochaines conÂf\'erences.

Au cours de cet expos\'e je n'ai pas l'intention de faire des maÂth\'ematiques, je ne ferai aucune d\'emonstration, mais simplement de parler sur les math\'ematiques. Or, il y a beaucoup de fa\c{c}ons d'envisager les math\'eÂmatiques, en raison des diff\'erentes \'ecoles de la philosophie des sciences, et \'egalement des sous-vari\'et\'es propres \`a la philosophie des math\'ematiques (logicisme, formalisme, intuitionnisme, platonisme, etc...). D\`es lors il est bon que je vous pr\'evienne que ce que je vais dire a plut\^ot tendance \`a orienter les choses vers le point de vue formaliste. Gardez seulement \`a l'esprit que ce point de vue est loin de mettre un point final \`a la quesÂtion {\it Que sont les math\'ematiques?}\\

\noindent{\large\bf I.	Les math\'ematiques naturelles}

Dans l'antiquit\'e et jusque vers la fin du 18e si\`ecle, les math\'ematiques \'etaient consid\'er\'ees comme l'expression des lois de la nature. La chose \'etait d'ailleurs si \'evidente que les auteurs n'abordent pas cette question en elle-m\^eme mais plut\^ot (comme Aristote ou Galil\'ee) la question de savoir si tous les ph\'enom\`enes naturels rel\`event des math\'eÂmatiques ou une partie seulement.
Le sens des symboles math\'ematiques \'etait unique et clair. Pascal \'ecrira vers 1660 [qu'on ne puisse tout d\'efinir et tout prouver] "c'est ce que la g\'eom\'etrie enseigne parfaitement. Elle ne d\'efinit aucune de ces choses espace, temps, mouvement, nombre, \'egalit\'e, ni les semÂblables qui sont en grand nombre, parce que ces termes-l\`a d\'esignent si naturellement les choses qu'ils signifient, \`a ceux qui entendent la langue, que l'\'eclair\-cissement qu'on voudrait faire apporterait plus d'obscurit\'e que d'instruction. Car il n'y a rien de plus faible que le discours de ceux qui veulent d\'efinir ces mots primitifs".
Ainsi, si Pascal consid\`ere qu'il est mauvais de chercher \`a d\'efinir les notions primitives des math\'ematiques, c'est simplement parce qu'il esÂtime que ces notions sont claires. Il n'envisage pas que le sens de ces notions primitives puisse n'avoir que peu d'importance pour les math\'ematiques elles-m\^emes.
Le premier penseur qui s'est rendu compte que le fait que les math\'ematiques nous apprenaient quelque chose sur la nature faisait probl\`eme si on le rapprochait du caract\`ere rigoureux du raisonnement math\'ematique, est incontestablement Kant. La solution qu'il propose pour r\'esoudre cette difficult\'e \`a savoir l'introduction de la c\'el\`ebre cat\'egorie des {\it jugements synth\'etiques a priori}, nous appara\^{\i}t aujourd'hui un peu comme une clause {\it ad hoc}, mais c'est qu'il \'etait impossible d'imaginer \`a cette \'epoque une math\'ematique qui ne soit pas naturelle.
Depuis l'antiquit\'e la g\'eom\'etrie \'etait la branche ma\^{\i}tresse des math\'ematiques et c'est pr\'ecis\'ement dans le domaine de la g\'eom\'etrie que vont na\^{\i}tre de nouvelles id\'ees qui vont se d\'evelopper tout au long du 19e si\`ecle et aboutir \`a une modification radicale de la conception des math\'ematiques.
C'est seulement \`a partir du 19e si\`ecle, que l'on trouve, en r\'eaction \`a ces nouvelles id\'ees, des math\'ematiciens partisans des math\'ematiques naturelles.
Si Hermite \'ecrit vers 1870 "Je crois que les nombres et les fonctions de l'Analyse ne sont pas le produit arbitraire de notre esprit, je pense qu'ils existent en dehors de nous avec le m\^eme caract\`ere de n\'ecessit\'e que les choses de la r\'ealit\'e objective, et nous les rencontrons et les \'etudions comme les physiciens, les chimistes et les zoologistes", c'est que ceci n'est plus compl\`etement \'evident. Que s'est-il produit ? L'apparition des g\'eom\'etries non-euclidiennes.\\

\noindent{\bf\large II.	Naissance du formalisme}

Dans les El\'ements d'Euclide (3e si\`ecle avant J.C.) l'axiome des parall\`eles est ainsi formul\'e "Si une ligne droite qui en rencontre deux autres forme d'un m\^eme c\^ot\'e avec les droites des angles internes dont la somme est moindre que deux droits, les deux premi\`eres droites se rencontrent ou leurs prolonge\-ments du c\^ot\'e o\`u la somme est inf\'erieure \`a deux droits".
Apr\`es de nombreuses tentatives souvent ing\'enieuses durant des si\`ecles pour d\'emontrer ce postulat \`a partir des autres axiomes, au d\'ebut du 19e si\`ecle quelques math\'ematiciens (Bolyai, Gauss, Lobachevski) en \'etaient arriv\'es \`a la conviction que les autres axiomes ne tranchaient pas entre cette hypoth\`ese et son contraire et certains d'entre eux affirmaient m\^eme qu'on pouvait \`a partir d'autres postulats et d\'evelopper tout un ensemble de cons\'equences "Une \'etrange g\'eom\'etrie, tout \`a fait diff\'erente de la n\^otre, enti\`erement coh\'erente en elle-m\^eme" dira Gauss vers 1820-1830.
Les g\'eom\'etries non-euclidiennes \'etaient n\'ees, on ne pouvait plus consid\'erer la g\'eom\'etrie classique que comme une parmi d'autres possibles et Riemann montrera en 1854 que certaines de ces diff\'erentes g\'eom\'etries planes peuvent se repr\'esenter dans l'espace comme la g\'eom\'etrie des g\'eod\'esiques d'une surÂface, de courbure totale positive ou n\'egative, d\'emontrant ainsi par un raisonnement de th\'eorie des mod\`eles avant la lettre, que ces nouvelles g\'eom\'etries sont aussi fiables que l'ancienne.

Le d\'eveloppement de la g\'eom\'etrie va s'intensifier dans le courant du 19e si\`ecle (Cayley, Pl\"{u}cker, Klein et le programme d'Erlangen, Pash, etc...). Progressivement \'emerge l'id\'ee que la validit\'e des raisonnements ne tient pas \`a la clart\'e du sens des notions primitives mais aux relations alg\'ebriques qu'elles entretiennent les unes avec les autres.
Cayley \'ecrira en 1859 "1) Le mot point peut signifier point et le mot droite peut siÂgnifier droite, 2) le mot point peut signifier droite et le mot droite signifier point. Nous pouvons par une telle extension comprendre les th\'eor\`emes corr\'elatifs sous un \'enonc\'e commun. Nous allons dans la g\'eom\'etrie \`a deux dimensions inclure la g\'eom\'etrie sph\'erique, les mots plan, point et droite signifiant \`a cette fin surface sph\'erique, arc (d'un grand cercle) et point (c'est-\`a-dire paire de points oppos\'es) de cette surface sph\'erique".
On traitait dans cet esprit d'un m\^eme coup un th\'eor\`eme sur les coniques et le th\'eor\`eme dual obtenu par transformation par polaires r\'eciproques.

Pash \'ecrira tr\`es clairement en 1882 dans ses le\c{c}ons sur la g\'eom\'etrie nouvelle "Il faut en effet, pour que la g\'eom\'etrie devienne une science d\'eductive, que la mani\`ere dont on tire les cons\'equences soit partout ind\'e\-pendante du sens des concepts g\'eom\'e\-triques, comme elle doit l'\^etre des figures, seuls sont \`a prendre en consid\'eration les rapports entre les concepts g\'eom\'etriques pos\'es par les propositions et les d\'efinitions emÂploy\'ees".
Cantor pourra dire en 1883 "La math\'ematique est enti\`erement libre dans son d\'eveloppement et ses concepts ne sont li\'es que par la n\'ecessit\'e d'\^etre non contradictoires et coordonn\'es aux concepts ant\'erieurement introduits par des d\'efinitions pr\'ecises".

Les choses vont assez rapidement \`a partir de 1880, une floraison de g\'eom\'etries appara\^{\i}t et un point culminant est atteint par Hilbert en 1899 dans ses {\it Fondements de la g\'eom\'etrie} o\`u il \'etudiera quasiment toutes les g\'eom\'etries possibles et montrera l'ind\'ependance des divers axiomes par le moyen d'interpr\'etations. La rigueur de ce travail et la puissance de sa g\'en\'eralit\'e (g\'eom\'etries non archim\'ediennes, non arguesiennes, etc.) a fortement impressionn\'e les contemporains, il apparaissait difficile d'aller plus loin.

Mais d\'ej\`a la m\'ethode axiomatique envahissait d'autres branches des math\'e\-matiques. En 1889, Peano donnait une version de ses c\'el\`ebres neuf axiomes de l'arithm\'etique. Quelques ann\'ees auparavant G. Frege proÂposait dans son {\it Begriffschrift} (1879) le premier langage compl\`etement formalis\'e (dans lequel sont introduits pour la premi\`ere fois les quantiÂficateurs existentiels et universels) et s'en servait dans ses {\it Grundsetze der Arithmetik }(1883) pour fonder l'arithm\'etique sur la logique pure.
Apr\`es quelques soubresauts dus \`a la d\'ecouverte par Russell en 1901 d'une contradiction (l'antinomie de Russell) dans un syst\`eme d'une aussi grande g\'en\'eralit\'e que celui de Frege, Russell et Whitehead propoÂsaient en 1910 dans leurs {\it Principia mathe\-matica} un syst\`eme formel qui se mettait \`a l'abri des paradoxes par l'\'elimination des d\'efinitions non pr\'edicatives.

Alors qu'au 19e si\`ecle les math\'ematiques se pr\'esentaient comme plusieurs domaines ind\'ependants g\'eom\'etrie, alg\`ebre, th\'eorie des foncÂtions dont le nombre s'\'etait accru r\'ecemment par l'apparition des g\'eoÂm\'etries non-euclidiennes, de la g\'eom\'etrie projective et de la th\'eorie des groupes, apr\`es la publication des Principia Mathematica l'ensemble des th\'eories math\'ematiques connues, depuis la g\'eom\'etrie jusqu'aux nombres transfinis de Cantor, est compl\`etement unifi\'ee en une seule th\'eorie, dont une forme plus commode utilis\'ee aujourd'hui est l'axiomatique de la th\'eoÂrie des ensembles de Zermelo-Fraenkel (ZF).\\

\noindent{\large\bf III.	Le formalisme de la th\'eorie des Ensembles}

Je vais maintenant pr\'eciser un peu le formalisme de la th\'eorie des ensembles, afin de tenter de vous convaincre en quelques minutes de deux choses pour lesquelles on estime g\'en\'eralement qu'elles ne peuvent vraiment s'acqu\'erir qu'en au moins un an de pratique personnelle.
La premi\`ere que j'exposerai sommairement est que le langage de la th\'eorie des ensembles est purement formel, on peut m\^eme dire mat\'eriel assemblages form\'es de signes parfaitement codifi\'es.
La seconde, que je ne ferai qu'\'evoquer, est que toute pens\'ee math\'ematique du moins la quasi-totalit\'e peut s'exprimer dans ce langage.

La th\'eorie des ensembles de Zermelo-Fraenkel s'exprime dans le langage de la logique des pr\'edicats du 1er ordre en posant un certain nombre d'axiomes sp\'ecifiques.
La logique des pr\'edicats du 1er ordre est une th\'eorie logique dans laquelle les quantificateurs $\forall\;\exists$ ne peuvent porter que sur des variables d'individus et non des variables de pr\'edicats (c'est-\`a-dire de propri\'et\'es).

\noindent{\bf Langage du 1er ordre}

	a)	Symboles de variables	$x,	y, z, x', y' , z' , x'', y'', z''\ldots$

	b)	Symboles de fonctions

		 - fonctions $0$-aires ou constantes,

		 - fonctions unaires	$f(\,), g(\,)\ldots$,

		 - fonctions binaires $h(\,,\,)\ldots$,

		 - fonctions n-aires	 $h(\,,\dots,\,)\ldots$,    

	c)	Symboles de pr\'edicats

		 - pr\'edicats unaires	$P(\,),\;Q(\,)\ldots$,

		 - pr\'edicats n-aires	$R(\,,\ldots,\,)\ldots$,

	d)	Les symboles $\neg$ (non), $\vee$ (ou), $\exists$ (il existe).

	e)	Un pr\'edicat binaire particulier not\'e $=$.

On d\'efinit alors les assemblages de signes qu'on appellera termes et ceux qu'on appellera formules.

\noindent{\bf Termes}	

i) une variable est un terme,

	ii)	si $u_1,\ldots,u_n$, sont des termes et si $f$ est une fonction n-aire
	alors $f(u_1,\ldots,u_n)$ est un terme,

\noindent{\bf Formules}

	a)	une formule atomique est un assemblage de la forme $P(a_1,\ldots,a_n)$	o\`u les $a_i$ sont des termes et $P$ un pr\'edicat n-aire.

b)	i) une formule atomique est une formule

\hspace{.5cm} ii)	si $\varphi$ est une formule, $\neg\varphi$ est une formule

\hspace{.5cm} iii)	si $\varphi$ et $\psi$ sont des formules, $\varphi\vee\psi$ est une formule,

\hspace{.5cm} iv)	si $\varphi$ est une formule, $\exists x\varphi$ est une formule et de m\^eme en rempla\c{c}ant $x$ par une autre variable.

Abr\'eviations on note $A\Rightarrow B$ pour  $\neg A\vee B$, $A\&B$  pour $\neg(A\Rightarrow\neg B)$, et enfin $\forall x A$ pour $\neg\exists x\neg A$.\\

\noindent{\large\bf Logique des pr\'edicats du 1er ordre}

On appelle ainsi une th\'eorie formalis\'ee utilisant le langage du 1er ordre plus un certain nombre d'axiomes et de r\`egles d'inf\'erence, perÂmettant de construire des formules \`a partir des axiomes, formules qu'on appellera des \underline{th\'eor\`emes}.\\

\noindent{\bf Axiomes}

Les formules suivantes sont des axiomes

1)	toute formule de la forme $A\vee\neg A$ o\`u $A$ est une formule.

2)	Si $A(x)$ est une formule ayant $x$ comme variable libre (non dans le champ d'un quantificateur) et si $a$ est un terme
$$A(a)\Rightarrow\exists x A(x)$$ est un axiome.

	3)	Si	$x$ est une variable
		 $x=x$ est un axiome,

	4)	Si $x_1,\ldots,x_n,y_1,\ldots, y_n$ sont des variables et $f$ une fonction n-aire
	$$(x_1=y_1)\&(x_2=y_2)\&\cdots\&(x_n=y_n)\Rightarrow f(x_1,\ldots,x_n)=f(y_1,\ldots, y_n)$$	est un axiome.

5)	et de m\^eme si $P$ est un pr\'edicat n-aire
	$$(x_1=y_1)\&(x_2=y_2)\&\cdots\&(x_n=y_n)\Rightarrow (P(x_1,\ldots,x_n)\Rightarrow P(y_1,\ldots, y_n))$$
	est un axiome.\\

\noindent{\bf R\`egles d'inf\'erence}

1)	de $A$ on d\'eduit $B\vee A$

	2)	de $A\vee A$	on d\'eduit $A$

3)	de $(A\vee B)\vee C$ on d\'eduit $A\vee(B\vee C)$

4)	de $A\vee B$ et de $\neg A\vee C$ on d\'eduit $B\vee C$

5)	Si $x$ n'est pas une variable libre de $B$, 

de $A\Rightarrow B$ on d\'eduit
$(\exists x A)\Rightarrow B$.\\

\noindent{\bf Th\'eorie des Ensembles}

Cette th\'eorie utilise un langage du 1er ordre sans symbole de fonction.

Un seul pr\'edicat binaire (en plus de =) not\'e $\in$ et un certain nombre d'axiomes pour le d\'etail desquels je renvoie \`a Krivine [10].

1)	Extensionalit\'e :
$$\forall z(z\in x\Leftrightarrow z\in y)\Rightarrow x=y$$

2)	R\'egularit\'e (ou fondation)
$$\exists y(y\in x)\Rightarrow \exists y(y\in x\; \&\;\neg\exists z(z\in x\&z\in y))$$

3)	Axiomes de partie : quelque soit la formule $A$
$$\exists z\forall x(x\in z\Leftrightarrow x\in y\;\&A)$$
est un axiome d\`es que les variables $x,y,z$ sont distinctes et
n'interviennent pas dans A.

4)	Le sch\'ema de remplacement.

5)	L'axiome de l'ensemble des parties.

6)	L'axiome de l'infini.

7)	(\'eventuellement) L'axiome du choix.\\

Pour ce qui est de voir comment on peut introduire \`a partir de cela toutes les math\'ematiques, nombres entiers, r\'eels, espaces topologiques, etc... je vous renvoie au cours d'analyse de C. WAGSCHAL, en atÂtirant seulement votre attention sur le fait que ce cours est r\'edig\'e dans le cadre de l'axiomatique de Zermelo qui est plus faible que ZF (le sch\'ema de rempla\-cement y est remplac\'e par l'axiome de compr\'ehension), mais d\'ej\`a suffisante pour faire la majeure partie de l'analyse.\\

\noindent{\large\bf IV.	Formalisation et ind\'ecidabilit\'e}

Ainsi nous sommes dans la situation suivante, l'ensemble des math\'ematiques se pr\'esente sous la forme d'un syst\`eme formel o\`u chaque th\'eor\`eme n'est qu'une cons\'equence purement logique et \'etrang\`ere \`a toute intuition sensible de certains axiomes. La v\'erit\'e de tel \'enonc\'e n'\'etant que relative \`a ce syst\`eme d'axiome, les math\'ematiques d'un c\^ot\'e se sont unifi\'ees, de l'autre ont perdu tout caract\`ere absolu les th\'eories ne contiennent que ce qu'on a mis dedans.
Si c'\'etait la seule conclusion \`a laquelle avait men\'e la formaÂlisation des math\'ematiques cela ne serait pas bien int\'eressant. Mais ceci a conduit les math\'ematiciens \`a consiÂd\'erer et \`a tenter d'\'etudier ces sortes de jeux que sont les syst\`emes formels eux-m\^emes, et ont d\'ecouvert des ph\'enom\`enes curieux auxquels en v\'erit\'e ils ne s'attendaient pas. Ce sera l'objet principal de tout ce cycle de conf\'erences.
Sans trop d\'eflorer les sujets des autres conf\'erenciers, je dirai simplement ceci\\

\noindent{\large\bf A.	Syntaxe}

A partir du moment o\`u toutes les math\'ematiques avaient \'et\'e ramen\'ees \`a des encha\^{\i}nements de signes, il ne devait pas \^etre bien sorcier de voir si ces encha\^{\i}nements, ind\'ependamment de toute signification atÂtribu\'ee aux symboles, peuvent mener a une contradiction, c'est-\`a-dire \`a $A\&\neg A$ ou non. C'est du moins ce qu'a cru Hilbert pendant un temps.
Cependant certains syst\`emes formels en particulier ceux obtenus a partir de la formalisation des math\'ematiques, quoique form\'es d'encha\^{\i}nements codifi\'es de suites de signes qu'on pourrait si l'on veut r\'ealiser mat\'eriellement en bois par exemple ou faire manipuler par un ordinateur, sont d'une telle richesse combinatoire, que sur certaines de leurs proÂpri\'et\'es plane un doute d'une nature particuli\`ere. La  th\`ese suivante sera pr\'ecis\'ee et comment\'ee par les prochaines conf\'erences
"Il y a des syst\`emes combinatoires parfaitement d\'etermin\'es (le hasard en est compl\`etement \'etranger) dont certaines propri\'et\'es (comme de savoir s'ils sont contradictoires ou non) sont \`a jamais en dehors du champ de la connaissance certaine en tout cas telle qu'elle \'etait entendue au 19e si\`ecle."
Cet aspect syntaxique sera illustr\'e et d\'evelopp\'e dans les conÂf\'erences de J.Y. GIRARD.\\

\noindent{\large\bf B.	S\'emantique}

La th\'eorie des ensembles de Zermelo-Fraenkel est un syst\`eme formel parti\-culier qui a l'avantage de constituer un cadre conceptuel tr\`es large dans lequel les probl\`emes de signification et de mod\`eles des syst\`emes formels peuvent \^etre abord\'es.
Un mod\`ele d'une th\'eorie formelle doit \^etre compris comme "une r\'ealit\'e" dont la th\'eorie serait la formalisation, c'est-\`a-dire comme un univers d'objets et de relations satisfaisant aux lois \'edict\'ees par les th\'eor\`emes de cette th\'eorie.
Un des r\'esultats fondamentaux de la th\'eorie des mod\`eles est que les syst\`emes formels ont le plus souvent plusieurs mod\`eles on dit qu'ils sont non cat\'egoriques.

Nous utilisons ces syst\`emes en attribuant un sens (le sens usuel) aux symboles, mais en fait d'autres significations sont possibles, d'autres objets que ceux dont on a voulu formaliser les propri\'et\'es saÂtisfont les \'enonc\'es.
Par exemple quoique l'on puisse faire toute la th\'eorie classique des nombres r\'eels dans Zermelo-Fraenkel, \`a partir de tout mod\`ele de Zermelo\-Fraenkel on peut construire un mod\`ele d\'enombrable dans lequel \'evidemment la relation d'appartenance et la propri\'et\'e d'\^etre d\'enombrable n'ont pas le sens habituel.

De m\^eme il existe des mod\`eles de l'arithm\'etique de Peano dans lesquels il y a des entiers infinis, et ces mod\`eles non-standards trouÂvent  des appli\-cations int\'eressantes.
La th\'eorie des mod\`eles a d\'evelopp\'e des concepts et des m\'ethodes par lesquels nos connaissances de la notion de symbole s'est trouv\'ee grandement affin\'ee par rapport \`a la r\'eflexion sp\'eculative de l'approche intuitive.

On peut dire si l'on veut que le math\'ematicien traditionnel se pr\'eoccupe de savoir s'il n'a pas d\'esign\'e le m\^eme objet par deux symboles diff\'erents (c'est ce qu'il fait lorsqu'il pose ou r\'esout une \'equation) alors que le logicien se pose la question de savoir si d'aventure un m\^eme symbole ne d\'esignerait pas plusieurs choses diff\'erentes.

Cet aspect s\'emantique sera illustr\'e particuli\`erement par l'exÂpos\'e de A. LOUVEAU sur l'hypoth\`ese du continu.

\newpage

\noindent{\LARGE\bf Le th\'eor\`eme d'incompl\'etude de G\"{o}del}

\noindent Jean-Yves Girard\\

Le formalisme, c'est-\`a-dire la position philosophique sur les math\'ematiques dont le repr\'esentant le plus \'eminent fut Hilbert, constitue le sujet essentiel de cet expos\'e. Plus pr\'ecis\'ement, nous allons exposer en d\'etail l'id\'eologie formaliste, dont la quintessence est le "programme de Hilbert" (en \'evoquant en contrepoint, une id\'eologie concurrente, l'intuitionnisme). Le programme de Hilbert avait ceci de particulier qu'il s'agissait d'un programme math\'e\-matique tr\`es pr\'ecis, tendant \`a d\'emontrer la justesse de l'ontologie formaliste (et non pas une vague succession d'intentions) c'\'etait s\^urement assez pr\'esomp\-tueux de la part de Hilbert que de vouloir trancher une question philosophique par le biais des maÂth\'ematiques, mais sans doute moins ridicule que l'inverse, \`a savoir d\'ecider une question scientifique au moyen de consid\'erations "philo\-sophiques" g\'en\'erales. En tout cas, on s'expose \`a la possibilit\'e de r\'efutation, car il est bien connu que les conjectures math\'ematiques peuvent \^etre d\'emontr\'ees (c'est ce qu'esp\'erait Hilbert pour son programme) ou r\'efut\'ees c'est ce qu'il advint. En 1931, par son th\'eor\`eme d'incompl\'etude, G\"{o}del r\'efutait clairement le programme de Hilbert. Cet expos\'e s'arr\^etera donc avec le th\'eor\`eme d'incompl\'etude. La semaine prochaine, nous nous pr\'eoccuperons de ce qu'il advint apr\`es car si le th\'eor\`eme d'incompl\'etude a montr\'e la fausset\'e des vues hilbertiennes sur les math\'ematiques, il ne nous a pas propos\'e grand-chose \`a la place, il y a eu des travaux dans la lign\'ee de Hilbert, et dont l'\'evaluation exacte reste d\'elicate, je pense en particulier aux travaux de Gentzen.\\

\noindent{\large\bf1.	La philosophie formaliste}

Les succ\`es de la formalisation (voir l'expos\'e pr\'ec\'edent) ont conduit \`a la position philosophique suivante, qui fut celle de Hilbert, mais qui a toujours de nombreux d\'efenseurs (notamment en France, voir les {\it El\'ements de Math\'ematique} de Bourbaki)
Les math\'ematiques doivent \^etre analys\'ees comme une activit\'e d\'epourvue de signification, comme, disons, le jeu d'\'echecs. Il s'agit, au moyen de r\`egles formelles fix\'ees \`a l'avance, de construire certains assemblages de symboles (les \'enonc\'es math\'ematiques et leurs d\'emons\-trations).

Faut-il en conclure que le formalisme sous-entend une attitude ontologique (l'ontologie est la partie de la philosophie qui traite de l'existence) coh\'erente ?  il n'en est rien, et nous pouvons facilement distinguer entre plusieurs positions formalistes

1 -	D'abord des gens pour qui rien n'existe. Cette position me semble difficilement conciliable avec l'id\'ee d'analyser les math\'ematiques comme un jeu formel, il semble en effet implicite dans cette id\'ee que les symboles formels existent vraiment on est donc amen\'e \`a reconna\^{\i}tre une certaine existence aux assemblages de symboles ainsi obtenus, et une certaine signi\-fication aux propri\'et\'es combinatoires (coh\'erence, etc...) de ces assemblages.

2 -	A l'oppos\'e se situe le monsieur qui dit au fond, toutes ces questions d'existence me passent au-dessus de la t\^ete, je ne veux pas savoir si les objets math\'ematiques existent ou non en tout cas, ce qui est s\^ur, c'est le c\^ot\'e purement formel mentionn\'e plus haut. Il ne s'agit pas vraiment d'une prise de position, mais plut\^ot d'une fuite...

3 -	La position de Hilbert, elle, est nettement plus int\'eressante\\

\noindent Les seuls objets qui existent pour Hilbert, ce sont des objets finitistes (ou encore \'el\'ementaires) il s'agit l\`a des \^etres r\'eels dont parlent les math\'ematiques. Quels sont-ils au juste ? Essentiellement les entiers, mais aussi les assemblages formels (qui peuvent \^etre repr\'esent\'es par des entiers).

Parmi les propri\'et\'es math\'ematiques, il en est qui n'ont pas de sens (comme par exemple des propri\'et\'es qui feraient r\'ef\'erence \`a des objets non finitistes), et d'autres qui ont une signification (et que l'on appellera encore "r\'eelles", "\'el\'ementaires", "finitistes"). Ces \'enonc\'es pleins de sens, Hilbert propose de les identifier aux propri\'et\'es qui ont lieu pour toute valeur des variables.

Exemples
	$$\forall x (x^2+2x+1=(x+1)^2)$$
$$\forall x\forall y\forall z\forall n(xyz\neq 0\;\&\;n>2\Longrightarrow x^n+y^n\neq z^n)$$
la conjecture de Riemann (dans ce cas, elle n'est pas directement \'el\'emen\-taire, mais on peut la ramener \`a une identit\'e, d'ailleurs sans int\'er\^et apparent). Parmi ces \'enonc\'es, il y a les formules de coh\'erence Coh(ZF) exprime que la th\'eorie ZF (Zermelo-Fraenkel) est ¥coh\'erente, c'est-\`a-dire

$$\forall\pi\;(\pi\;{\mbox{d\'emonstration dans ZF}} \Rightarrow\pi\;{\mbox{ne se termine pas par}}\;0=1)$$
enfin, ce ne sont pas toutes les d\'emonstrations math\'ematiques (y compris quand elles ne font intervenir que des \'enonc\'es et des objets \'el\'ementaires) qui seront irr\'eprochables, il faut que les principes utilis\'es dans les d\'emonstrations soient particuli\`erement imm\'ediats, on parlera alors de d\'emons\-tration \'el\'ementaire (ou finitiste).

Les math\'ematiques (et sp\'ecialement celles du XXe si\`ecle) sont loin de se r\'eduire \`a ce noyau d\'elimit\'e par Hilbert, nous avons des objets abstraits (ou id\'eaux) (pensons aux ultrafiltres, aux espaces de... Hilbert). Pour Hilbert, ils n'existent pas, de m\^eme les propri\'et\'es abstraites (ou id\'eales) n'ont pas de sens, par exemple, la n\'egation d'une identit\'e n'a pas, pour Hilbert, de signification, de m\^eme quant aux principes abstraits qui se retrouvent dans de nombreuses d\'emonstrations, par exemple, l'axiome du choix. Cette distinction entre "r\'eel" et "abstrait" va nous conduire au programme de Hilbert.\\

\noindent{\large\bf 2.	Le programme de Hilbert}

L'id\'ee de Hilbert est de montrer qu'il est possible, en th\'eoÂrie du moins, de se passer des objets abstraits, des \'enonc\'es abstraits. Voil\`a donc le programme qu'il se proposait de r\'ealiser.

Soit $R$ une propri\'et\'e \'el\'ementaire, \'etablie \`a l'aide de m\'ethodes abstraites, montrer que $R$ peut \^etre \'etablie \`a l'aide de m\'ethodes \'el\'eÂmentaires.

Ce programme, d'\'enonc\'e simple, appelle un certain nombre de remarques

- D'abord, il se place dans une certaine tradition, rappelons-nous le th\'eor\`eme des nombres premiers (Hadamard, de la Vall\'ee-Poussin), d\'emontr\'e au moyen de la th\'eorie des fonctions analytiques, mais qui a \'et\'e \'el\'ementaris\'e au XXe si\`ecle, on a trouv\'e des d\'emonstrations plus proches du r\'esultat, c'est-\`a-dire n'utilisant pas de principes abstraits comme des r\'esultats g\'en\'eraux sur des fonctions analytiques... Le principe de "puret\'e des m\'ethodes" dit qu'il est bien, \'el\'egant... de se restreindre, lors d'une d\'emonstration, \`a des principes et des \'enonc\'es en rapport avec la concluÂsion, une d\'emonstration abstraite d'un \'enonce r\'eel viole la puret\'e des m\'ethodes, et le programme de Hilbert nous donne donc un principe de puret\'e des m\'ethodes g\'en\'eralis\'e.

- Le programme de Hilbert suppose la compl\'etude des math\'ematiques \'el\'e\-men\-taires, en effet, il ne fait pas de doute que toute identit\'e vraie est d\'emontrable par des m\'ethodes abstraites ad hoc, et donc par le proÂgramme de Hilbert, par des m\'ethodes \'el\'ementaires. D'ailleurs, les formaÂlistes attard\'es utilisent "vrai" pour "d\'emontrable" et "faux" pour "r\'efutable". Les ph\'enom\`enes de compl\'etude (c'est-\`a-dire l'ad\'equation entre
vrai et\break "prou\-vable") ne sont pas si rares que cela, mentionnons quelques exemples

-	la th\'eorie des corps r\'eels clos est compl\`ete,

-	la compl\'etude est vraie pour les \'enonces co-\'el\'ementaires (n\'egations d'\'enonc\'es \'el\'ementaires) en fait tout \'enonc\'e co-\'el\'ementaire vrai est prouvable dans les math\'ematiques "co-\'el\'ementaires"... malheureusement, ce n'est pas \`a cette classe-l\`a que Hilbert s'est int\'eress\'e pour son programme.

Une d\'emonstration valable du programme de Hilbert ne saurait \^etre qu'\'el\'e\-mentaire~: sinon, on ne pourrait rien conclure, puisque ce qui n'est pas \'el\'ementaire n'a pas de sens... Concr\`etement, la d\'emonstration du programme de Hilbert consiste \`a se placer dans une th\'eorie abstraite (par exemple ZF), et d'y consid\'erer un \'enonc\'e \'el\'ementaire, et une d\'emonsÂtration abstraite T de cet \'enonc\'e R, il s'agirait alors, au moyen de transformations \'el\'emen\-taires sur T, d'obtenir une d\'emonstration \'el\'eÂmentaire T' du m\^eme \'enonc\'e R. Ce genre de r\'esultat (dans un autre contexte) n'a rien d'utopique on sait par exemple (\`a l'aide des r\'esultats de G\"{o}del sur l'hypoth\`ese du continu) transformer, au moyen d'un proc\'ed\'e
tout ce qu'il y a d'\'el\'ementaire, toute d\'emonstration dans ZF+AC+GCH d'un \'enonc\'e \'el\'ementaire, en une d\'emonstration dans ZF (mais la d\'emonsÂtration n'est pas \'el\'ementaire, et l'analogie s'arr\^ete l\`a...). L'\'ecole de Hilbert se fit d'abord les dents sur de tout petits syst\`emes (fragments de l'arithm\'etique), en attendant d'\^etre \`a m\^eme de d\'emon\-trer le r\'esultat pour la th\'eorie des ensembles. Mais G\"{o}del donna un coup d'arr\^et brutal au programme...\\

\noindent{\bf\large3.	D\'emonstration de coh\'erence}

Une formulation \'equivalente du programme de Hilbert, est la recherche de d\'emonstrations de coh\'erence \'el\'ementaires~:

(i)	si le programme de Hilbert est r\'ealis\'e, supposons que par exemple, ZF soit contradictoire, c'est-\`a-dire d\'emontre O = 1 , qui est un \'enonc\'e \'el\'eÂmentaire, par le programme de Hilbert, cet \'enonc\'e peut donc \^etre d\'emontr\'e par des m\'ethodes \'el\'ementaires mais ces m\'ethodes ne sont dites "\'el\'ementaires" que parce que nous n'avons aucun doute en ce qui les concerne, en particulier elles ne sauraient mener \`a une absurdit\'e du genre de O = 1. Nous venons ainsi d'obtenir une d\'emonstration \'el\'ementaire de coh\'erence pour ZF (la d\'emonstration \'el\'ementaire du programme de Hilbert pour ZF, et l'\'evidence \'el\'ementaire de la coh\'erence des math\'ematiques \'el\'ementaires (si la coh\'erence des math\'ematiques \'el\'ementaires n'est pas \'etablie de mani\`ere \'el\'ementaire, ici encore, une application du programme de Hilbert rendra notre d\'emons\-tration \'el\'ementaire) )

ii)	r\'eciproquement, si la th\'eorie ZF est coh\'erente, prenons un \'enonc\'e \'el\'ementaire, disons le th\'eor\`eme de Fermat, et supposons l'y avoir prouv\'e;
si par ailleurs, les entiers $a,b,c, \neq 0$ et $n>2$ sont tels que $a^n+b^n=c^n$ par exemple si $17^{101}+23^{101} = 24^{101}$ , alors l'\'equation $a^n+b^n=c^n$ serait prouvable en th\'eorie des ensembles ZF, ce qui donnerait une contradiction dans ZF. Nous venons de donner informellement une d\'emonstration \'el\'ementaire du programme de Hilbert (dans le cas du th\'eor\`eme de Fermat), en donnant un proc\'ed\'e \'el\'ementaire pour transformer une d\'emonstration non-\'el\'ementaire de Fermat dans ZF en une d\'emonstration \'el\'ementaire, au moyen d'une d\'emons\-tration \'el\'ementaire de coh\'erence de ZF...\\

\noindent{\large\bf 4. L'intuitionnisme}

On comprendra mieux le formalisme en parlant de son fr\`ere ennemi, l'intui\-tionnisme, dont le champion fut Brouwer. Les deux th\'eories se reÂjoignaient dans un m\^eme refus du r\'ealisme (c'est-\`a-dire de la philosophie "na\"{\i}ve" des math\'e\-matiques, pour laquelle les \'enonc\'es math\'e\-matiques parlent d'objets r\'eels dans un univers id\'eal \`a pr\'eciser...). Par contre, les deux th\'eories divergeaient radicalement pour ce qui est de l'analyse de l'acÂtivit\'e math\'ematique.

-	Pour Hilbert, l'activit\'e math\'ematique est m\'ecanique, pour lui le math\'e\-maticien id\'eal serait un robot, essayant syst\'ematiquement toutes les d\'emonstrations \'el\'ementaires possibles. Tout ce que nous appelons "intuition" n'est rien d'autre que ce qui permet au math\'ematicien de se hausser au niveau de ce math\'ematicien id\'eal...

-	au contraire, pour Brouwer, l'activit\'e math\'ematique est un acte cr\'eateur de l'esprit, qui ne saurait en aucune fa\c{c}on \^etre m\'ecanis\'e, en particulier, pour lui, les math\'ematiques \'etaient par essence non formalisables.

Le syst\`eme intuitionniste se d\'eveloppe par des principes logiÂques totalement diff\'erents des principes classiques, notamment le refus du principe du tiers exclu, tout cela se justifie, suivant Brouwer, en faisant r\'ef\'erence \`a un "sujet pensant" qui produit des r\'esultats math\'ematiques. Une discussion approfondie de l'intuitionnisme serait ici hors de propos contentons-nous de remarquer que cette th\'eorie, bien qu'exÂtravagante par certains c\^ot\'es, pr\'esente une vision de l'activit\'e math\'emaÂtique moins d\'esolante que celle de Hilbert, malheureusement, l'intuitionnisme n'a pas tenu ses promesses...\\

\noindent{\large\bf 5.	R\'efutation du programme de Hilbert (1931)}

Le premier th\'eor\`eme d'incompl\'etude de G\"{o}del associe \`a toute th\'eoÂrie coh\'erente T satisfaisant des conditions tr\`es g\'en\'erales, un \'enonc\'e \'el\'eÂmentaire C, ayant les propri\'et\'es suivantes

- C est vrai

-	C n'est pas d\'emontrable dans T.

Ceci d\'etruit le programme de Hilbert, car si on prend pour th\'eoÂrie T une th\'eorie coh\'erente contenant toutes les math\'ematiques \'el\'ementaires,
(par exemple, on pourrait prendre pour T l'arithm\'etique de Peano AP), l'\'eÂnonc\'e de G\"{o}del C est 
\'el\'ementaire, d\'emontr\'e par le th\'eor\`eme de G\"{o}del (puisque le th\'eor\`eme dit qu'il est vrai) et la coh\'erence de T, mais sans d\'emonstration dans T, donc sans d\'emonstration \'el\'ementaire.

On trouve plus souvent la r\'ef\'erence au deuxi\`eme th\'eor\`eme d'incompl\'etude, qui dit la chose suivante : sous des conditions presque aussi g\'en\'erales que celles du premier th\'eor\`eme, l'\'enonc\'e \'el\'ementaire qui exprime la coh\'erence de T n'est pas d\'emontrable dans T. Ce deuxi\`eme th\'eor\`eme d\'etruit le programme sous sa forme "coh\'erence" mais il n'est pas n\'ecessaire (sa d\'emonstration est plus difficile que celle du premier) \`a la r\'efutation du programme de Hilbert le premier th\'eor\`eme suffit.

Le programme de Hilbert est donc r\'efut\'e en 1931, mais, par un ph\'enom\`ene psychologique bien compr\'ehensible, le formalisme a gard\'e une grande partie de ses adeptes : en effet le formalisme pr\'esente une vision assez simple des math\'ematiques (le m\'ecanisme, la compl\'etude, les objets abstraits comme des fa\c{c}ons de parler, mais n'existant pas...), et on sait bien que les visions simplistes, r\'eductrices, du monde ont toujours un impact sans commune mesure avec leur succ\`es r\'eel. Tout aurait \'et\'e bien diff\'erent si au lieu de d\'emolir le programme de Hilbert sans rien mettre \`a sa place qui ait les m\^emes vertus d'attraction, G\"{o}del avait trouv\'e une caract\'erisation simple de la prouvabilit\'e en termes de v\'erit\'e, ou le conÂtraire... Nous allons maintenant voir d'un peu plus pr\`es ces th\'eor\`emes d'incompl\'etude.\\

\noindent{\large\bf6.	Le premier th\'eor\`eme d'incompl\'etude}

Nous allons nous int\'eresser ici uniquement \`a l'arithm\'etique de Peano AP. La premi\`ere chose est de repr\'esenter les donn\'ees syntaxiques de AP au moyen d'entiers. Cela ne pr\'esente th\'eoriquement aucune difficult\'e, vu que le langage, et les constructions qui gravitent autour, forment des ensembles d\'enombrables. La bijection que l'on \'etablit entre \'enonc\'es et $\mathbb{N}$ s'appelle num\'erotation de G\"{o}del. Rappelons que le langage de l'arithm\'e\-tique est bas\'e sur $0,S,+,.,=,<$ \`a chacun de ces symboles, ainsi qu'aux connecÂteurs, aux quantificateurs et aux variables, on attribue un entier son
\underline{nombre de} \underline{G\"{o}del}~: 
$\ulcorner 0\urcorner=1\;;\;\ulcorner  S\urcorner=3\;;\;\ulcorner +\urcorner=5\;;\; \ulcorner.\urcorner=7\;;\;\ulcorner=\urcorner=9\;;\;\ulcorner<\urcorner=11\;;\;\ulcorner\neg\urcorner=13\;;\;\ulcorner\vee\urcorner=15\;;\;\ulcorner\&\urcorner=17\;;\; \ulcorner\Rightarrow\urcorner=19\;;\; \ulcorner\exists\urcorner=21\;;\; \ulcorner\forall\urcorner=23\;;\; \ulcorner x_n\urcorner=25+2n.$

\noindent Les nombres de G\"{o}del ne sont en rien remarquables ; on pourrait imaginer beaucoup d'autres mani\`eres de num\'eroter les symboles du langage. Il ne s'agit pas d'entiers qui auraient des relations cach\'ees avec les symboles qu'ils repr\'esentent, je dis cela \`a cause de l'association qu'on peut \^etre amen\'e \`a faire avec d'autres "nombres deÉ" (Bernoulli par exemple). Si $a_0,\ldots,a_{n-1}$ est une suite d'entiers de longueur n, on introduit un entier
$\ulcorner a_0,\ldots,a_{n-1}\urcorner=p_1^{a_0+1}\cdots p_n^{a_{n-1}+1}$ o\`u $p_n$ d\'esigne le n-ime nombre premier, alors $\ulcorner a_0,\ldots,a_{n-1}\urcorner$ d\'etermine n
uniquement, ainsi que les entiers $a_0,\ldots, a_{n-1}$.

On sait qu'il est possible, en utilisant la "notation polonaise", c'est\-\`a-dire en \'ecrivant $\vee AB$ \`a la place de $A\vee B$, $\&AB$ \`a la place de $A \& B$,... d'\'eliminer les parenth\`eses, et tout \'enonc\'e du langage appara\^{\i}t comme une suite finie de symboles $0,S,\ldots,x_n,\ldots$ dans un certain ordre. Si $A$ est un tel \'enonc\'e,
on peut donc \'ecrire $A =s_0\ldots s_m$ o\`u $s_0\ldots s_m$ sont des symboles  $0,S,\ldots,x_n,\ldots$
et on posera donc $\ulcorner A\urcorner = \ulcorner \ulcorner s_0\urcorner,\ldots,\ulcorner s_{m}\urcorner\urcorner$. Une d\'emonstration de A dans AP, c'est une suite $A_0,\ldots,A_p$ d'\'enonc\'es, telle que $A_p = A$, et la suite est
b\^atie en respectant les r\`egles d'inf\'erence du calcul des pr\'edicats, et les
axiomes de Peano (dont l'\'enonc\'e pr\'ecis ne joue aucun r\^ole ici). On peut donc associer a une d\'emonstration $A_0,\ldots,A_p$ un nombre de G\"{o}del $\ulcorner A\urcorner = \ulcorner \ulcorner A_0\urcorner,\ldots,\ulcorner A_{p}\urcorner\urcorner$.

On consid\`ere la propri\'et\'e $Dem(d,a)$ : $d$ est le nombre de G\"{o}del d'une d\'emonstration de l'\'enonc\'e dont le nombre de G\"{o}del est $a$. Nous allons \'enoncer quelques \'evidences (techniquement, cela demande des notions de r\'ecursivit\'e, et ce n'est pas du tout \'evident)

- $Dem (d,a)$ est une formule de l'arithm\'etique (on voit bien que les r\`egles syntaxiques de formation des \'enonc\'es et des d\'emonstrations vont \^etre expressibles, au moyen du nombre de G\"{o}del, par des propri\'et\'es arithm\'eÂtiques)

- surtout d\`es que $Dem(d,a)$ est vrai, il est d\'emontrable dans AP, cela peut se voir ainsi : une d\'emonstration est une suite finie de symboles, qui v\'erifie un certain nombre de crit\`eres combinatoires que l'on peut v\'eÂrifier m\'ecaniquement ; or tout ce qui est m\'ecanique est du ressort de la d\'emons\-trabilit\'e (je n'ai pas dit que l'on peut d\'ecider m\'ecaniquement la d\'emons\-trabilit\'e dans AP), c'est-\`a-dire que le processus m\'ecanique de v\'eriÂfication qu'un assemblage de symboles est une d\'emonstration peut \^etre confi\'e
\`a une th\'eorie formelle comme AP. En fait les \'enonc\'es, qui comme Dem(d,a) sont prouvables d\`es qu'ils sont vrais, sont les co-\'el\'ementaires. (On dit plus souvent $\Sigma_1^0$)

L'\'etape suivante consiste \`a reproduire le fameux \underline{paradoxe du menteur}, connu depuis l'Antiquit\'e : "je mens" ; \`a proprement parler, le paradoxe du menteur donnerait "je ne suis pas vrai" ; ici nous allons conÂsid\'erer "je ne suis pas prouvable", \`a savoir un \'enonc\'e $A$ \'equivalent \`a $\forall d\neg Dem(d,\ulcorner A\urcorner)$ en fait on a mieux que l'\'equivalence $A$ et  $\forall d\neg Dem(d,\ulcorner A\urcorner)$ sont le m\^eme \'enonc\'e. Un tel \'enonc\'e $A$ se construit facilement avec un peu
plus de technique que ce que nous avons introduit... Remarquons que $A$ est
par construction un \'enonc\'e \'el\'ementaire (on dit encore $\Pi_1^0$ ).

$A$ n'est pas prouvable dans AP : en effet, si $A$ \'etait prouvable dans AP, on y prouverait $A$ et donc $Dem(d,\ulcorner A\urcorner)$ pour un certain $d$ ; mais comme $A$ est $\forall d\neg Dem(d,\ulcorner A\urcorner)$, on obtiendrait une contradiction dans AP donc si AP est coh\'erente, l'\'enonc\'e $A$ n'est pas d\'emontrable dans AP. Mais comme $A$ est \'equivalent \`a sa non-prouvabilit\'e dans AP, et que nous venons de voir que $A$ n'est pas d\'emontrable dans AP, il s'ensuit que $A$ est vrai, mais non d\'emontrable. Remarquons que rien n'exclut que $\neg A$ soit d\'emontrable dans AP (sinon notre intuition que les th\'eor\`emes de AP sont vrais). Par une am\'elioration due \`a Rosser (1936), on peut construire $A$ vrai et \'el\'eÂmentaire tel que ni $A$ ni $\neg  A$ ne soient d\'emontrables dans AP.

Une autre remarque est le r\'esultat suivant, d\^u \`a Tarski, on ne peut pas trouver de pr\'edicat de v\'erit\'e dans AP, \`a savoir d'\'enonc\'e $V$ tel que $V(\ulcorner B\urcorner) \Leftrightarrow (B$ soit d\'emontrable dans AP) pour tout $B$ il suffit de refaire la construction de G\"{o}del en rempla\c{c}ant $\exists d\;Dem(d,\ulcorner B\urcorner)$ par $V(\ulcorner B\urcorner)$.\\

\noindent{\bf\large7.	Deuxi\`eme th\'eor\`eme d'incompl\'etude}

$Coh({\tt AP})$ est l'\'enonc\'e $\forall d\neg Dem(d,\ulcorner0=S0\urcorner)$, qui exprime que AP est coh\'e\-rente~; il est extr\^emement p\'enible de d\'emontrer le premier th\'eor\`eme d'incom\-pl\'etude formellement dans AP, mais on y arrive : si AP est coh\'erente, $A$ est vrai~: l'\'enonc\'e $Coh({\tt AP})\Rightarrow A$ est donc d\'emontrable dans AP, et comme $A$ n'est pas d\'emontrable dans AP, $Coh({\tt AP})$ ne peut pas \^etre d\'emontrable dans AP, si AP est coh\'erente. L'id\'ee est simple, mais la mise en oeuvre math\'e\-matique est tr\`es lourde...\\

\noindent{\large\bf8. Signification des th\'eor\`emes d'incompl\'etude}

Le r\'esultat de G\"{o}del est extr\^emement c\'el\`ebre, d'ailleurs surtout hors des milieux math\'ematiques. Nous avons d\'ej\`a eu l'occasion de dire que ce th\'eor\`eme n'avait pas affect\'e de fa\c{c}on sensible la conception des math\'e\-matiques du math\'ematicien moyen. La renomm\'ee de ce r\'esultat dans des cercles non math\'e\-matiques vient surtout d'une m\'ecompr\'ehension : on retient l'imiÂtation du paradoxe du menteur, et cela donne quelque chose comme (si on parle du deuxi\`eme th\'eor\`eme, qui a de loin le plus de succ\`es) "une th\'eorie ne peut pas se penser elle-m\^eme", ou pire, en extrapolant "la science ne peut pas se penser elle-m\^eme". Ce qui est certain, c'est qu'un r\'esultat comme le th\'eor\`eme de G\"{o}del appelle naturellement de nombreux commentaires et des extrapolations vari\'ees : c'est qu'on a l'impression que son imporÂtance d\'epasse de loin le strict cadre technique o\`u il se place... Ceci dit, ce n'est pas une raison pour dire n'importe quoi.

- Quand on dit qu'une th\'eorie ne peut pas se penser elle-m\^eme on fait un contresens sur la signification du th\'eor\`eme, les \'enonc\'es qui expriment dans AP les propri\'et\'es de AP, telle la formule $Coh({\tt AP})$, montrent bien que AP peut "se penser elle-m\^eme" ; si on veut dire par l\`a qu'il y a des \'enonc\'es sur AP qui ne sont pas d\'emontrables dans AP elle-m\^eme, il conÂviendra de remarquer qu'il y a aussi des \'enonc\'es vrais qui ne parlent pas de AP, et qui ne sont pas d\'emontrables dans AP, ainsi que des propri\'et\'es vraies de AP qui sont d\'emontrables dans AP : l'accent mis sur la coh\'erence est arbi\-traire...

- Quand l'extrapolation va plus loin que les th\'eories formelles, \`a savoir quand on veut faire dire au th\'eor\`eme d'incompl\'etude quelque chose sur la pens\'ee en g\'en\'eral, on passe loin de la siÂgnification profonde du r\'esultat, en effet le th\'eor\`eme d'incompl\'etude est avant tout un r\'esultat sur la pens\'ee m\'ecanique (une th\'eorie formelle) et s'il semble honn\^ete d'extrapoler les r\'esultats d'incompl\'etude au domaine de la pens\'ee m\'ecanique (ordinateurs, jeux...), appliquer le th\'eor\`eme d'incompl\'etude \`a l'activit\'e humaine est hautement farfelu, d'ailleurs, quand je dis appliquer, c'est un bien grand mot : on cite G\"{o}del \`a la place o\`u nagu\`ere on aurait cit\'e, disons... Mallarm\'e, dont les vers \'enigmatiques renferment, on le sait, tout ce qu'il faut pour ce genre de situation : il importe seulement de ne pas m\'elanger les genres et de ne pas utiliser les
r\'esultats scientifiques sous une forme tellement m\'etaphoriques qu'ils en sont devenus m\'econnaissables.

Pour r\'esumer, si nous devons dire en une phrase quelle est la signification de la r\'efutation par G\"{o}del du programme de Hilbert, c'est~:

"il y a des choses qui ne sont pas du ressort du m\'ecanisme." 

Voil\`a ce que disent les th\'eor\`emes d'incompl\'etude, et ce message est r\'ep\'et\'e sous de nombreuses formes techniques, je pense aux th\'eor\`emes d'ind\'eci\-dabilit\'e (voir 5\`eme conf\'erence). S'il y a des choses qui ne sont pas du ressort du m\'ecanisme, cela veut dire qu'il y a quelque chose de non m\'ecanique dans, disons, la pens\'ee math\'ematique par exemple, la production de nouveaux axiomes est n\'ecessairement chaotique (vue du point de vue d'un robot) il s'agit de parier \`a chaque \'etape sur l'\'etape suivante : il n'est pas question de produire des axiomes par une m\'ethode r\'eguli\`ere.\\

\noindent{\large\bf9.	D\'emonstrations de coh\'erence relative}(\footnote{Ce qui suit constitue la 2e conf\'erence de J.Y. GIRARD})

Apr\`es le r\'esultat de G\"{o}del, il devenait impossible de penser \`a raisonner sur la coh\'erence d'une mani\`ere qui soit \'epist\'emologiquement absolument convaincante. Comment penser un instant qu'une d\'emonstration de coh\'erence, qui utilise plus que l'arithm\'etique, a une r\'eelle valeur, au niveau de la philosophie des math\'ematiques... En fait ceci ne concerne que les d\'emons\-trations de coh\'erence, dans ce qu'elles ont de plus pr\'etenÂtieux : les d\'emonstrations de coh\'erence absolue. Par contre, il est une autre cat\'egorie de d\'emonstrations de coh\'erence, les d\'emons\-trations de coh\'erence relative, c'est-\`a-dire des \'enonc\'es de la forme
$$Coh({\tt T}) \Rightarrow Coh({\tt U})$$
(si T est coh\'erente, alors U l'est). La valeur \'epist\'emologique de tels r\'esultats est plus limit\'ee, mais par contre, vu que de telles d\'emonstraÂtions se font toujours dans les math\'ematiques \'el\'ementaires ch\`eres \`a Hilbert, cette valeur \'epist\'emologique est incontestable ; par exemple, le r\'esultat de 1939 de G\"{o}del : $Coh({\tt ZF}) \Rightarrow Coh({\tt ZF + AC + GCH})$
(coh\'erence relative de l'axiome du choix et de l'hypoth\`ese du continu
par rapport \`a ZF) montre que l'axiome du choix et l'hypoth\`ese du continu ne sont pas plus "risqu\'ees" que la th\'eorie des ensembles~; nous avons ici une vraie r\'eduction... De la m\^eme mani\`ere, les r\'esultats de Cohen (coh\'eÂrence relative de la n\'egation de l'axiome du choix, de la n\'egation de l'hypoth\`ese du continu), montrent que l'on peut aussi faire des choix oppos\'es qui ne sont pas plus risqu\'es. Mais ne d\'eflorons pas trop le sujet de la conf\'erence d'A. LOUVEAU.\\

\noindent{\large\bf10.	Signification des d\'emonstrations de coh\'erence}

Un mot de Kreisel "les doutes quant \`a la coh\'erence sont plus douteux que la coh\'erence elle-m\^eme". Au fond, personne ne doute vraiment de la coh\'erence de l'arithm\'etique ; s'il existe un probl\`eme li\'e aux th\'eoÂr\`emes de G\"{o}del, il importe de ne pas dramatiser \`a l'exc\`es, et de ne pas appliquer \`a ces situations d\'elicates des consid\'erations trop brutales.

Un des r\'esultats les plus remarquables de Gentzen est son introÂduction, dans les ann\'ees 30 du calcul des s\'equents, et la d\'emonstration pour le calcul des s\'equents, du \underline{Hauptsatz} (qui est un v\'eritable principe de puret\'e des m\'ethodes pour le calcul des pr\'edicats, c'est-\`a-dire la logique pure) n'est pas sans cons\'equences pour les questions de coh\'erence. Avant de parler du calcul des s\'equents, expliquons comment on peut par exemple utiliser le Hauptsatz pour \'etablir, en un sens presque satisfaisant, la coh\'erence des math\'ematiques \'el\'ementaires.

(ARP) l'arithm\'etique r\'ecursive primitive est une bonne formulation des math\'eÂmatiques \'el\'ementaires ; ce syst\`eme est form\'e ainsi

- il y a un grand nombre (une infinit\'e) de lettres de fonctions, pour repr\'esenter des fonctions combinatoires \'el\'ementaires (fonctions r\'ecursives primitives) comme $+, . ,\exp$, etc..., et pour chacune de ces fonctions,
on \'ecrit des axiomes de d\'efinition, par exemple

\noindent$2^0=1$

\noindent$2^{Sx}=2^x+2^x$

- on permet l'induction (r\'ecurrence) sur les \'enonc\'es de la forme $t = u$, o\`u $t$ et $u$ sont des \underline{termes}, c'est-\`a-dire sont construits en utilisant
les lettres de fonction mentionn\'ees plus haut, et des variables.

Ce syst\`eme est une bonne formalisation des math\'ematiques \'el\'ementaires~; mais est-il coh\'erent ? Il y a une d\'emonstration "na\"{\i}ve" de la coh\'erence de ARP, qui consiste seulement \`a remarquer que tout th\'eor\`eme de ARP est vrai, par exemple, \`a montrer que si $\forall x\exists y \;t[x,y]=0$ est d\'emontrable dans ARP, alors, pour tout $x$, on peut trouver un $y$ tel que $t[x,y]$ soit \'egal \`a $0$. Malheureusement, on voit bien que dans la d\'efinition de la v\'erit\'e, on utilise des quantificateurs qui nous font sortir de la classe des \'enonc\'es \'el\'ementaires typiquement, le "on peut trouver y" est un quantificateur ind\'esirable... Heureusement, le Hauptsatz peut ici \^etre utilis\'e, on remarque d'abord que tous les axiomes de ARP sont de la forme $t = u$, y compris l'axiome (les axiomes en fait) d'induction. Le Hauptsatz nous dit alors (c'est un principe de puret\'e des m\'ethodes) qu'une d\'emonstration de $0 = S0$ se place dans la partie purement \'equationnelle de ARP, c'est-\`a-dire qu'\`a partir des \'equations qui nous servent d'axiomes, utilisant la transitivit\'e et la sym\'etrie de l'\'egalit\'e,  on a d\'eduit $0=S0$ !. Cette r\'eduction (qui utilise le Hauptsatz) de ARP \`a sa partie \'equationnelle est faite par des moyens purement \'el\'ementaires. Par contre, il faut d\'emontrer la coh\'erence de la partie purement \'equationnelle (sans logique) de ARP, et pour cela, nous avons une fonction de valuation, qui \`a tout terme clos (sans variable) de ARP, associe sa valeur~; bien entendu, cette fonction de valuation ne saurait \^etre parmi les fonctions de ARP. Cette partie de la d\'emonstration de coh\'erence de ARP n'est donc pas, conform\'ement au th\'eor\`eme de G\"{o}del, formalisable dans ARP. La fin de la d\'emonstration est simple, on d\'emontre que, si $t$ et $u$ sont des termes clos tels que l'\'equation $t = u$ soit d\'emontrable dans la partie \'equationnelle de ARP, alors $V(t)=V(u)$; pour les axiomes c'est imm\'ediat, et la transitivit\'e, la sym\'etrie de l'\'egaÂlit\'e, pr\'eservent cette propri\'et\'e... En particulier, comme $V(0) = 0 \neq 1 =V(S0)$, $0=S0$ n'est pas d\'emontrable dans le calcul \'equationnel, et donc pas d\'emontrable tout court dans ARP.

Nous voyons ici que la partie probl\'ematique de la coh\'erence (tout ce qui touche aux quantificateurs) est d\'emontr\'ee de mani\`ere \'el\'eÂmentaire, tandis que le th\'eor\`eme de G\"{o}del ne s'applique vraiment qu'\`a la partie purement \'equationnelle de ARP, pour laquelle nous n'avons pas de doute s\'erieux.\\

\noindent{\large\bf11.	Le calcul des s\'equents}

Il est difficile ici de ne pas \'evoquer la m\'emoire de Jacques Herbrand, mort pr\'ematur\'ement en 1931, et dont les travaux (le th\'eor\`eme de Herbrand) pr\'efigurent les r\'esultats de Gentzen...

Pour obtenir son principe de puret\'e des m\'ethodes, Gentzen modiÂfie les notions de base de la logique : la notion essentielle n'est plus le concept d'\'enonc\'e, mais celui de s\'equent, Gentzen appelle s\'equent une expression formelle $\Gamma\truc\Delta$ , o\`u $\Gamma$ et $\Delta$ sont des suites finies d'\'enonÂc\'es. La signification intuitive de $A_1,\ldots,A_n\truc B_1,\ldots,B_m$	c'est que
si tous les $A_i$ sont vrais, alors un des $B_j$ l'est~: 
$$A_1\&\ldots\& A_n\Rightarrow B_1\vee\ldots B_m$$

(En particulier, $\truc A$ veut dire A, et $A \truc$ veut dire $\neg A$, et puis le s\'equent vide $\truc$ veut dire l'absurdit\'e, c'est-\`a-dire $0 = S0$). Les r\`egles du calcul des s\'equents sont tout \`a fait remarquables ; on les divise en quatre groupes
\begin{center}
	(I) AXIOMES
\end{center}
$$A\truc A\quad(A\quad \mbox{\'enonc\'e quelconque})$$
\begin{center}
	(Il) REGLES STRUCTURELLES
\end{center}

\noindent Affaiblissement
$$\frac{\Gamma\truc\Delta}{\Gamma\truc A,\Delta}\quad dA\qquad\qquad\frac{\Gamma\truc\Delta}{\Gamma,A\truc \Delta}\quad gA$$

\noindent Echange
	$$\frac{\Gamma\truc\Delta^\prime,A,B,\Delta^{\prime\prime}}{\Gamma\truc \Delta^\prime,B,A,\Delta^{\prime\prime}}\quad dE\qquad\qquad
\frac{\Gamma^\prime,A,B,\Gamma^{\prime\prime}\truc\Delta}{\Gamma^\prime,B,A,\Gamma^{\prime\prime}\truc \Delta}\quad gE$$

\noindent Contraction
	$$\frac{\Gamma\truc A,A,\Delta}{\Gamma\truc A,\Delta}\quad dC\qquad\qquad
\frac{\Gamma,A,A\truc\Delta}{\Gamma,A\truc \Delta}\quad gC$$
\begin{center}
(III) REGLES LOGIQUES
\end{center}

\noindent R\`egles du \&
	$$\frac{\Gamma\truc A,\Delta\quad \Gamma^\prime\truc B,\Delta^\prime}
{\Gamma,\Gamma^\prime\truc A\&B,\Delta,\Delta^\prime}\;\;d\&\qquad\frac{\Gamma,A\truc\Delta}
{\Gamma,A\&B\truc\Delta}\;g1\&\qquad\frac{\Gamma,B\truc \Delta}
{\Gamma,A\&B\truc\Delta}\;g2\&$$

\noindent R\`egles du $\vee$
	$$\frac{\Gamma\truc A,\Delta}{\Gamma\truc A\vee B,\Delta}\;d1\vee\qquad
\frac{\Gamma\truc B,\Delta}{\Gamma\truc A\vee B, \Delta}\;d2\vee\qquad
\frac{\Gamma,A\truc \Delta\;\quad\Gamma^\prime,B\truc\Delta^\prime}
{\Gamma,\Gamma^\prime,A\vee B\truc\Delta,\Delta^\prime}\;g\vee$$
\noindent R\`egles du $\Rightarrow$
	$$\frac{\Gamma,A\truc B,\Delta}{\Gamma\truc A\Rightarrow B,\Delta}\;d\Rightarrow\qquad
\frac{\Gamma\truc A,\Delta\;\quad\Gamma^\prime,B\truc \Delta^\prime}
{\Gamma,\Gamma^\prime,A\Rightarrow B\truc \Delta,\Delta^\prime}\;g\Rightarrow$$
\noindent R\`egles du $\neg$
	$$\frac{\Gamma,A\truc\Delta}{\Gamma\truc \neg A, \Delta}\;d\neg\qquad\qquad
\frac{\Gamma\truc A,\Delta}{\Gamma,\neg A\truc\Delta}\;g\neg$$
\noindent R\`egles du $\forall$
	$$\frac{\Gamma\truc A,\Delta}{\Gamma\truc \forall xA,\Delta}\;d\forall(^\ast)\qquad\qquad
\frac{\Gamma,A[t]\truc\Delta}{\Gamma,\forall xA\truc \Delta}\;g\forall$$
\noindent R\`egles du $\exists$
$$\frac{\Gamma\truc A[t],\Delta}{\Gamma\truc \exists xA,\Delta}\;d\exists\qquad
\frac{\Gamma,A\truc \Delta}{\Gamma,\exists xA\truc \Delta}\;g\exists(^\ast)$$
($^\ast$) r\`egles restreintes au cas o\`u $x$ n'est pas libre dans $\Gamma\truc \Delta$.
	\begin{center}
(IV) COUPURE
\end{center}
	$$\frac{\Gamma\truc A,\Delta\qquad\Gamma^\prime,A\truc \delta^\prime}
{\Gamma,\Gamma^\prime\truc \Delta,\Delta^\prime}\;{\tt C}$$
une d\'emonstration dans le calcul des s\'equents LK, c'est une suite de s\'e\-quents (dispos\'ee en forme d'arbre) construite \`a partir des axiomes, au moyen des r\`egles (II), (III), (IV) . Un mot sur les r\`egles

(II) les r\`egles structurelles sont des r\`egles assez anodines : elles permettent de faire quelques manipulations combinatoires simples (permu\-tation, augmentation, identification d'\'enonc\'es du m\^eme c\^ot\'e du signe $\truc$) ; il importe cependant de ne pas m\'epriser le pouvoir de la r\`egle de contraction, 

(III) les r\`egles logiques d\'emarquent la gen\`ese des \'enonc\'es ; on parle encore de r\`egles g\'en\'etiques pour la partie (I) (II) (III) de LK (calcul sans coupures), pour cette raison. C'est une r\'eussite ind\'eniable de Gentzen que d'avoir isol\'e, au moyen des s\'equents de tels principes "purs" de d\'emons\-tration.

(IV) la r\`egle de coupure exprime quant \`a elle la \underline{transitivit\'e} de la notion de cons\'equence logique ; c'est une r\`egle qui correspond \`a la pratique maÂth\'ematique chaque fois que dans un raisonnement math\'ematique, vous utilisez un r\'esultat g\'en\'eral, vous effectuez une coupure ; je donne un exemple simple~:

si on vous demande de calculer $499^2$ , il est peu vraisemblable que vous effectuiez le calcul il est tellement mieux de faire
$$499^2=500^2-2\times500+1$$

Pour cela, vous utilisez un r\'esultat g\'en\'eral $\forall x\;(x-1)^2=x^2-2x+1$ que vous appliquer au cas particulier de 499 ; ce processus se traduit par une coupure : soit $\pi$ la d\'emonstration classique de l'identit\'e susmen\-tionn\'ee,
alors, on peut obtenir $499 = 500 - 2.500 + 1$ ainsi
$$
\frac{
\begin{array}{c}
\pi\\
\vdots\\
\truc \forall x\;(x-1)^2=x^2-2x+1
\end{array}\qquad
\frac{(500-1)^2=500^2-2.500+1\truc (500-1)^2=500^2-2.500+1}
{\forall x\;(x-1)^2=x^2-2x+1\truc (500-1)^2=500^2-2.500+1}g\forall}
{\truc (500-1)^2=500^2-2.500+1}\;{\tt C}
$$
Les r\'esultats sur le calcul LK sont les suivants~:

1/ Un r\'esultat attendu, LK est une bonne formulation de la logique, autrement dit $A$ est d\'emontrable dans le calcul des pr\'edicats (encore : $A$ est vrai dans tout mod\`ele) si et seulement si $\truc A$ est d\'emontrable dans LK.

2/ Le Hauptsatz~: si un s\'equent est d\'emontrable dans LK il est aussi d\'eÂmontrable sans coupures, autrement dit la r\`egle de coupure est redondante.

Le corollaire le plus connu du Hauptsatz, c'est la \underline{propri\'et\'e} de la\break \underline{sous-formule}. Gentzen appelle sous-formule de $A $, $A$ lui-m\^eme et

-	si $A$ est $\neg B$, les sous-formules de $B$,

-	si $A$ est $B\&C$, $B\vee C$, $B\Rightarrow C$, les sous-formules de $B$ et/ou de $C$,

-	si $A$ est $\forall xB[x]$, $\exists xB[x]$, les sous-formules des $B[t]$, o\`u $t$ est un terme.
On voit qu'il y a toujours (d\`es que A contient un quantificateur), une infinit\'e de sous-formules distinctes pour A. On a la propri\'et\'e\\

Dans une d\'emonstration sans coupures de $\Gamma\truc \Delta$, tous les \'enonc\'es utilis\'es sont des sous-formules d'\'enonc\'es apparaissant dans $\Gamma$ ou dans $A$.\\

Ce r\'esultat est une cons\'equence des deux remarques suivantes

-	dans les r\`egles (II) et (III), tous les \'enonc\'es utilis\'es dans une quelconque pr\'emisse de la r\`egle, sont des sous-formules d'\'enonc\'es apparaisÂsant dans la conclusion,

-	la notion de sous-formule est transitive.

Ce principe, qui est un authentique principe de puret\'e des m\'eÂthodes pour le calcul des pr\'edicats (\`a un d\'etail pr\`es il y a des substiÂtutions arbitraires de termes pour fabriquer les sous-formules, mais on se rappellera que l'on peut se ramener \`a des calculs sans lettres de fonctions, o\`u les seuls termes sont donc des variables), nous a permis d'"\'eliminer les quanti\-ficateurs" dans ARP.

Le th\'eor\`eme de Gentzen est \'el\'ementaire, c'est-\`a-dire que Gentzen donne un proc\'ed\'e effectif pour \'eliminer les coupures. (Il est \`a remarquer que ce proc\'ed\'e, quoique effectif, d\'epasse de loin, les capacit\'es th\'eoriques
des ordinateurs !). En particulier, on peut se demander ce que donne une
d\'emonstration sans coupures du r\'esultat donn\'e plus haut (sur $499^2$) : c'est
tout simplement effectuer le calcul on voit que les d\'emonstrations sans coupures ont tendance \`a devenir longues et b\^etes, c'est la r\`egle de coupure qui concentre l'intelligibilit\'e des d\'emonstrations, et c'est pourquoi elle correspond \`a la pratique ; les r\`egles sans coupures, correspondent elles, \`a l'\'etude abstraite du raisonnement. Elles sont commodes \`a \'etudier, mais pas \`a manier.

Un mot sur la \underline{s\'emantique} des r\`egles sans coupures : c'est une s\'emantique \`a trois valeurs (vrai, faux, ind\'etermin\'e) ; les r\`egles (I) (Il) (III) sont valides pour l'interpr\'etation : si $A_1,\ldots,A_n$ sont vrais, alors 
un des \'enonc\'es $B_1,\ldots,B_m$ n'est pas faux. Pour ce type d'interpr\'etation,
la r\`egle de coupure n'est pas valide, puisque cela exigerait que $A$ non faux implique $A$ vrai, autrement dit la r\`egle de coupure caract\'erise la s\'emantique binaire, o\`u tout est vrai ou faux...

La post\'erit\'e de Gentzen est impressionnante en th\'eorie de la d\'emonstration les th\'eor\`emes d'\'elimination des coupures continuent \`a fournir une grande partie de l'inspiration et de la probl\'ematique du sujet...\\

\noindent{\large\bf12.	D\'emonstrations de coh\'erence pour l'arithm\'etique}

Gentzen a aussi donn\'e des d\'emonstrations de coh\'erence pour l'arithm\'etique~; ces d\'emonstrations reposent sur un principe additionnel, qui est un principe d'induction (r\'ecurrence) transfinie. Consid\'erons les
polyn\^omes exponentiels, construits, \`a partir de $0$, de la somme, et de
l'exponentiation $x^{P(x)}$ ; on peut comparer deux tels polyn\^omes par
$$P<Q\Leftrightarrow \exists n(\forall m>n\;(P(m)>Q(m))).$$
 On obtient un ordre lin\'eaire qui est un bon ordre (pas de suite infinie d\'ecroissante), et que l'on note $\varepsilon_0$. Gentzen d\'emontre la coh\'erence de AP (l'arithm\'etique de Peano) par induction transfinie jusqu'\`a  $\varepsilon_0$. L'id\'ee de la d\'emonstration est de r\'ep\'eter ce que
nous avons fait pour ARP, c'est-\`a-dire de se ramener \`a la partie purement
\'equationnelle de AP, qui est trivialement coh\'erente~; pour cela, il nous suffirait d'un principe de puret\'e des m\'ethodes, c'est-\`a-dire d'une g\'en\'eÂralisation du Hauptsatz \`a l'arithm\'etique. Mais le Hauptsatz n'est vrai que pour des syst\`emes sans axiomes (ou des axiomes de complexit\'e tr\`es simple comme ARP) ; l'id\'ee est de se ramener \`a cette situation en faisant appara\^{\i}tre le principe d'induction comme quelque chose de purement logique\footnote{Ceci est \`a rapprocher de la mani\`ere dont Bourbaki \'ecrit l'axiome du choix~: par un tour de passe-passe, cet axiome devient une cons\'equence purement logique des autres~; l'avantage est que cela permet (?) de couper court \`a des discussions sur le choix des axiomes qui obligeraient \`a mettre `math\'ematique' au pluriel...}

On \'ecrit donc les r\`egles
$$\frac{\cdots\Gamma\truc A[n],\Delta\cdots}
 {\Gamma\truc\forall xA[x],\Delta}\;d\forall\omega\qquad\qquad
 \frac{\Gamma,A[n_0]\truc \Delta}{\Gamma,\forall xA[x]\truc\Delta}\;g\forall\omega$$
et
$$\frac{\Gamma\truc A[n_0],\Delta}{\Gamma\truc \exists A[x],\Delta}\;d\exists\omega\qquad\qquad
 \frac{\cdots\Gamma,A[n]\truc \Delta\cdots}{\Gamma,\exists xA[x)\truc \Delta}\;g\exists\omega$$
\`a la place des r\`egles logiques \'ecrites plus haut pour $\forall$ et $\exists$ ; observer que $d\forall\omega$ et $g\exists\omega$ sont des r\`egles infinies~: les d\'emonstrations sont donc des arbres bien fond\'es, et on a affaire \`a une g\'en\'eralisation de la notion habituelle de d\'emonstration, qui n'a rien \`a voir avec la notion de d\'emons\-tration formelle. (Sch\"{u}tte, \`a qui on doit cette formulation, parle
de syst\`eme "semi-formel"). Rappelez-vous la justification du principe d'induction comme on a pu vous l'apprendre il y a quelques ann\'ees~: si $A[0]$ et $\forall x(A[x]\Rightarrow A[Sx])$, alors $A[0]\Rightarrow A[S0]$ , donc $A[S0]$ , puis $A[S0] \Rightarrow A[SS0]$ ; donc $A[SS0]$ etc..., autrement dit les s\'equents $A[0], \forall x(A[x]\Rightarrow A[Sx]) \truc A[n]$ sont prouvables pour tout n, et par la r\`egle $(d\forall\omega)$, $A[0],\forall x(A[x]\Rightarrow A[Sx])\truc \forall xA[x]$ est d\'emontrable ce qui donne bien une d\'emonstration purement logique du principe d'induction usuelle. Mais le calcul semi-formel mentionn\'e plus haut v\'erifie l'\'elimiÂnation des coupures, et de plus, si on part d'une d\'emonstration (finie) dans l'arithm\'etique, on obtient finalement une d\'emonstration sans coupures du m\^eme s\'equent, et de \underline{hauteur} born\'ee par $\varepsilon_0$ ; en particulier, en prenant le s\'equent $\truc$, on voit imm\'ediatement qu'une d\'emonstration sans coupures de ce s\'equent (qui se placerait dans la partie \'equationnelle de AP) ne
saurait exister ; mais ici, on utilise un principe ext\'erieur \`a AP		l'induction transfinie jusqu'\`a $\varepsilon_0$,	c'est le tribut qu'il faut payer au deuxi\`eme th\'eor\`eme d'incompl\'etude.\\

\noindent{\large\bf13.	Signification des d\'emonstrations de coh\'erence}

Un math\'ematicien fran\c{c}ais a dit : "Gentzen, c'est le type qui a d\'emontr\'e la coh\'erence de l'arithm\'etique, c'est-\`a-dire de l'induction transfinie jusqu'\`a $\omega$, au moyen de l'induction transfinie jusqu'\`a $\varepsilon_0$". Cette plaisanterie exprime tr\`es bien le malaise que l'on ressent devant le r\'eÂsultat de Gentzen~: qu'est-ce que ces r\'esultats veulent bien dire ?

1/ A y regarder de plus pr\`es, le travail de Gentzen n'est pas si ridicule que cela ; en effet, l'arithm\'etique, c'est l'induction jusqu'\`a $\omega$ sur des \'enonc\'es arbitraires, alors que dans le travail de Gentzen, comme l'a fait remarquer Kreisel, l'induction jusqu'\`a $\varepsilon_0$ n'est utilis\'ee que sur des \'enonc\'es sans quantificateurs, c'est-\`a-dire qu'il s'agit d'une induction \'el\'ementaire. Avec cette remarque, le travail de Gentzen peut tout \`a fait \^etre accept\'e comme une g\'en\'eralisation l\'egitime du programme de Hilbert, mais il n'en reste pas moins que la valeur \'epist\'emologique en reste limit\'ee, si on se borne \`a la question des d\'emonstrations de coh\'erence absolue ; apr\`es tout, l'induction jusqu'\`a $\omega$ sur des \'enonc\'es quelconques est plus pr\`es de notre intuition que l'induction (m\^eme \'el\'ementaire) jusqu'\`a $\varepsilon_0$ !

2/ Comme il n'est donc pas possible de consid\'erer les r\'esultats de Gentzen
comme des d\'emonstrations absolument irr\'efutables de la coh\'erence, Kreisel
a propos\'e de leur trouver un contenu positif (ind\'ependant de consid\'erations
oiseuses), c'est-\`a-dire des corollaires math\'ematiques indiscutables. Par
exemple des propri\'et\'es de puret\'e des m\'ethodes g\'en\'eralis\'ees, telles~:
 
Si $A$ est d\'emontrable dans AP, il est d\'emontrable en utilisant l'induction transfinie jusqu'\`a un ordinal $<\varepsilon_0$ et sur un \'enonc\'e de complexit\'e moindre que $A$ : ce r\'esultat est une application directe du r\'esultat d'\'elimination des coupures pour la $\omega$-logique.

3/ Peut-on quand m\^eme essayer de faire plus que de trouver de simples "contenus positifs" aux r\'esultats de Gentzen ? Il semble que oui : si $\alpha$ est
un ordinal, on peut d\'efinir\footnote{Ne pas confondre cette exponentielle ordinale avec l'exponentiation cardinale qui elle, 	est reli\'ee \`a la notation $A^B$ pour d\'esigner l'ensemble des applications de $B$ dans $A$.} l'exponentielle $2^\alpha$ , par~:

$2^0=1$

$2^{\alpha+1}=2^\alpha+2^\alpha$

$2^\lambda=\sup_{\mu<\lambda}2^\mu\qquad(\lambda\mbox{ limite })$

\noindent La "philosophie" du travail de Oentzen, c'est que

un quantificateur = une exponentielle

\noindent autrement dit on peut \'eliminer des quantificateurs au moyen d'exponentielles ordinales, ou le contraire. D'ailleurs, $\varepsilon_0$ n'est rien d'autre que le premier $\alpha\neq\omega$ tel que $\alpha=2^\alpha$, c'est-\`a-dire le r\'esultat de $\omega$ exponentiations...

Tout cela peut \^etre rendu plus explicite par le r\'esultat suivant~:
 les deux principes suivants sont \'equivalents dans les math\'ematiques \'el\'eÂmentaires
 
(1)	le fait que tout sous-ensemble de $\mathbb{N}^2$ a une projection sur $\mathbb{N}$.

(2)	le fait que si $\alpha$ est un bon ordre, $2^\alpha$ est un bon ordre.

L'int\'er\^et de cette \'equivalence est que, dans (1), l'op\'eration de
projection n'est pas \'el\'ementaire ; \'etant donn\'e $X\subset\mathbb{N}^2$ , comment d\'eterminer
les \'el\'ements de $pr_1(X)$ et de son compl\'ementaire sans utiliser d'op\'erations infinies~? Par contre, si on conna\^{\i}t $\alpha$, on conna\^{\i}t ipso facto $2^\alpha$ : ce qui fait probl\`eme, c'est alors le fait que l'exponentielle pr\'eserve le c\^ot\'e bien ordonn\'e des ordres lin\'eaires. On a donc une mise en relation entre deux
approches : d'une part une approche "r\'ealiste" o\`u il y a des objets, et surtout des constructions
 infinies, et qui correspond \`a notre pratique maÂth\'ematique, d'autre part, 
une approche "\'el\'ementaire", o\`u seules un certain type de construction est 
autoris\'e ; certaines propri\'et\'es non-\'el\'ementaires des cons\-tructions (ici \^etre
 bien ordonn\'e) permettent de rendre compte des op\'erations non \'el\'ementaires 
de l'univers "r\'ealiste".

Tout cela nous engage fortement \`a chercher une autre ontologie pour les 
math\'ematiques, c'est-\`a-dire \`a chercher les "vrais" objets de la praÂtique 
math\'ematique ailleurs que dans notre intuition premi\`ere. L'erreur de Hilbert 
semble \^etre d'avoir voulu faire des \'enonc\'es, des propri\'et\'es \'el\'eÂmentaires 
quelque chose de m\'ecanique, alors qu'il appara\^{\i}t clairement que de nouveaux 
principes \'el\'ementaires (comme l'induction jusqu'\`a $\varepsilon_0$ sur des \'equations) sont 
sans arr\^et n\'ecessit\'es pour rendre compte de la comÂplexit\'e logique croissante 
des math\'ematiques. Oui, l'univers math\'ematique se ram\`ene cer\-tainement \`a un 
sous-univers \'el\'ementaire ; mais ce dernier n'est pas plus simple que l'univers tout entier.
 
 \newpage
 \noindent{\LARGE\bf Ind\'ecidabilit\'e de l'hypoth\`ese du continu}
 
 \noindent Alain Louveau\\
 
\noindent{\large\bf 1.	L'hypoth\`ese du continu}

Dans ce troisi\`eme volet sur l'ind\'ecidabilit\'e en math\'ematiques, je vais parler de l'ind\'ecidabilit\'e en th\'eorie des ensembles. Le th\'eor\`eme d'incompl\'etude de G\"{o}del, appliqu\'e \`a cette th\'eorie, nous assure l'existence d'\'enonc\'es ni prouvables ni r\'efutables dans cette th\'eorie. Le second th\'eor\`eme donne un exemple de tel \'enonc\'e celui qui affirme la coh\'erence de la th\'eoÂrie elle-m\^eme. Ces r\'esultats sont certainement tr\`es importants du point de vue th\'eorique, mais d'un effet tr\`es faible sur la pratique math\'ematique. Tant que dans sa propre pratique, un math\'ematicien n'a pas rencontr\'e un \'enonc\'e qu'il cherchait \`a d\'emontrer et qui s'av\`ere ind\'ecidabie, les r\'esultats d'ind\'ecidabilit\'e restent sans effets. Mais les derni\`eres vingt ann\'ees ont vu fleurir, dans diverses branches des math\'ematiques, des r\'esultats d'ind\'epenÂdance. Et c'est de l'histoire de l'un d'entre eux - le plus c\'el\`ebre - que je vais parler~: l'hypoth\`ese du continu.

A la fin du 19e si\`ecle, le math\'ematicien Georg Cantor fonde la th\'eorie des ensembles par ses travaux sur la notion de cardinalit\'e~: deux ensembles ont la m\^eme cardinalit\'e si leurs \'el\'ements peuvent \^etre mis en correspondance biunivoque; C'est la notion naturelle de taille pour les ensembles finis. L'extension aux ensembles infinis a des propri\'et\'es un peu plus \'etranges par exemple l'ensemble des entiers, celui des entiers pairs et l'ensemble $\mathbb{Q}$ des rationnels ont la m\^eme cardinalit\'e (la plus petite pour les ensembles infinis, et qui est not\'ee $\aleph_0$ ). Un r\'esultat fondamental de Cantor assure qu'un ensemble $A$ n'a jamais la m\^eme cardinalit\'e que l'ensemble $\mathcal{P}(A)$ de toutes ses parties. En it\'erant l'op\'eration des parties, on construit donc des ensembles de cardinalit\'es de plus en plus grandes.
Ainsi la cardinalit\'e de $\mathcal{P}(\mathbb{N})$, not\'ee $2^{\aleph_0}$ , et qui est la m\^eme que celle de l'ensemble $\mathbb{R}$ des nombres r\'eels, est plus grande que $\aleph_0$ , etc... Il est naturel de se demander s'il y a des cardinalit\'es interm\'ediaires, en particulier s'il existe une partie de $\mathbb{R} $ n'ayant ni la cardinalit\'e $\aleph_0$ de $\mathbb{N}$, ni celle de $\mathbb{R}$. Cantor ne put construire de tel ensemble, et conjectura que cela devait \^etre impossible. C'est cet \'enonc\'e "Il n'y a pas de partie de $\mathbb{R}$
de cardinalit\'e interm\'ediaire entre$\aleph_0$ et $2^{\aleph_0}$, qui est connu sous le nom d'hypoth\`ese du continu (la droite r\'eelle, dans les textes anciens, est aussi appel\'ee continu).

Ce probl\`eme donna lieu \`a un grand nombre de recherches, avant
1938, ann\'ee o\`u K. G\"{o}del apporta une demi-r\'eponse~: il n'est pas possible
de r\'efuter l'hypoth\`ese du continu dans la th\'eorie ZFC\footnote{ZFC est la th\'eorie des ensembles de Zermelo-Fraenkel avec axiome du choix c'est la th\'eorie axiomatique usuelle de la pratique math\'ematique.}. Le r\'esultat
est un r\'esultat de coh\'erence relative~: s'il est possible de trouver une contra\-diction \`a partir de la th\'eorie des ensembles augment\'ee de l'hypoth\`ese du continu, alors la th\'eorie des ensembles, \`a elle seule, est d\'ej\`a contradictoire. L'hypoth\`ese du continu peut donc \^etre consid\'er\'ee comme un axiome suppl\'e\-mentaire "s\^ur". Mais est-ce vraiment un axiome suppl\'ementaire (et non pas un th\'eor\`eme)~? La r\'eponse attendit vingt-cinq ans. En 1963, P. Cohen prouva que pas plus qu'elle n'est r\'efutable, l'hypoth\`ese du continu n'est d\'emontrable nous sommes bien en pr\'esence d'un \'enonc\'e ind\'ecidable.

Les d\'emonstrations de G\"{o}del et Cohen sont n\'ecessairement techniques. Je vais cependant essayer d'en donner les id\'ees principales, en commen\c{c}ant par d\'egager le principe des d\'emonstrations d'ind\'ependance~: la construction de mod\`eles.\\

\noindent{\large\bf2.	Syntaxe et s\'emantique}

Pour illustrer ce qui va suivre, prenons une th\'eorie plus simple, par exemple la th\'eorie des corps. Le langage de cette th\'eorie comprend des symboles pour $+,\times, 0$ et $1$. Et les axiomes sont ceux des corps. Consid\'erons l'\'enonc\'e E de ce langage "$1 + 1\;=\;0$". Cet \'enonc\'e est-il d\'emontrable ou r\'efutable \`a partir des axiomes de la th\'eorie des corps ? Il s'agit d'un probl\`eme combinatoire portant sur des assemblages de symboles. Pourtant une premi\`ere r\'eponse, de nature tr\`es diff\'erente, vient \`a l'esprit. Si "$1 + 1\; =\;0$" \'etait d\'emontrable \`a partir des axiomes, il serait s\^urement vrai dans le corps $\mathbb{R}$, ce qui n'est pas le cas. E ne doit donc pas \^etre d\'emontrable. En disant cela, on est pass\'e de l'aspect syntaxique \`a l'aspect s\'emantique des th\'eories formelles. Etant donn\'ee une th\'eorie T dans
un certain langage, on appelle mod\`ele de T un ensemble, muni des relations et fonctions idoines, sp\'ecifi\'ees par le langage, et qui satisfait les axiomes de la th\'eorie T. Par exemple un mod\`ele de la th\'eorie des corps est un corps. Une fois que l'on poss\`ede cette notion s\'emantique, on introduit naturellement, \`a c\^ot\'e de la notion de cons\'equence syntaxique (existence d'une d\'emonstration), une notion de cons\'equence s\'emantique un \'enonc\'e du langage est cons\'equence s\'emantique des axiomes de T s'il est satisfait par tous les mod\`eles de T. Le fait que ces deux notions de cons\'equence co\"{\i}ncident est un autre remarquable r\'esultat de K. G\"{o}del, appel\'e le th\'eor\`eme de compl\'etude. Il a comme cons\'equence qu'une th\'eorie formelle est coh\'erente si et seulement si elle admet un mod\`ele. Par suite, pour revenir \`a notre \'enonc\'e E du langage de la th\'eorie des corps, il suffit, pour prouver l'ind\'ecidabilit\'e de E, de produire deux corps l'un, comme $\mathbb{R}$, dans lequel la n\'egation de E est satisfaite, et l'autre, par exemple $\mathbb{Z}/2$ dans lequel l'\'enonc\'e E est satisfait.

La technique utilis\'ee pour l'ind\'ependance de l'hypoth\`ese du continu va \^etre la m\^eme~: produire deux mod\`eles de ZFC, l'un satisfaisant l'hypoth\`ese du continu, l'autre sa n\'egation. Avec une difficult\'e, nous n'allons pas construire ces mod\`eles \`a partir de rien, car ce n'est pas possible en effet, ces mod\`eles seraient en particulier mod\`eles de la th\'eorie des ensembles, et par le th\'eor\`eme de compl\'etude la construction fournirait alors une preuve de la coh\'erence de cette th\'eorie, ce qui n'est pas possible par des arguments formalisables dans la th\'eorie des ensembles (ce sont les seuls que nous allons employer), par le th\'eor\`eme d'incompl\'etude. Ce n'est pas non plus ce que nous cherchons nous voulons \'etablir une coh\'erence relative \`a celle de la th\'eorie des ensembles. Nous allons donc partir d'un mod\`ele arbitraire de la th\'eorie des ensembles, suppos\'e exister, et construire \`a partir de lui les deux mod\`eles cherch\'es. Bien entendu, ces constructions d\'ependent d'une certaine connaissance g\'en\'erale des mod\`eles de la th\'eorie des ensembles.\\

\noindent{\large\bf3.	A quoi ressemble un univers ?}

La th\'eorie des ensembles va intervenir \`a deux niveaux dans ce qui suit d'une part les constructions effectu\'ees vont \^etre (au moins th\'eoriquement) formalisables dans cette th\'eorie. D'autre part, les objets
construits seront (on l'esp\`ere) des mod\`eles de la th\'eorie des ensembles.

Il faut donc bien s\'eparer ces deux niveaux, en particulier en ce qui concerne la terminologie. Dans le premier cas, nous parlerons d'ensembles intuitifs. Par exemple, un mod\`ele de la th\'eorie des ensembles, que nous appellerons un univers, est un ensemble intuitif, comme tout mod\`ele d'une th\'eorie. Et le langage de la th\'eorie des ensembles sp\'ecifiant une relation binaire, l'univers U est muni d'une telle relation, que nous noterons $\in^U$ peut donc \^etre repr\'esent\'e comme un graphe du type
suivant, o\`u 
\begin{center}
  \includegraphics[scale=0.6]{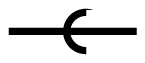}
\end{center}
  repr\'esente la relation $\in^U$.

  \begin{center}
    \includegraphics[scale=1.2]{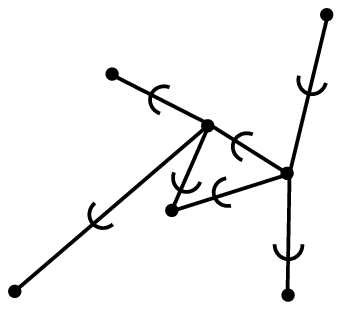}
  \end{center}
Nous appellerons ensembles les \'el\'ements (intuitifs) de U. Pour illustrer la distinction entre ensembles et ensembles intuitifs, consid\'erons le fragment suivant de l'univers U

\begin{center}
  \includegraphics[scale=1.2]{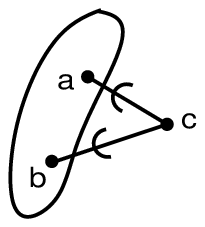}
\end{center}
Dans ce dessin, l'\'el\'ement $c$ de U est l'ensemble $\{a,b\}^{\tt U}$ . Il a
exactement pour \'el\'ements, au sens de U , les ensembles $a$ et $b$. Par
contre la "patate" dessin\'ee est la paire intuitive $\{a,b\}$. Il y a l\`a un ph\'enom\`ene g\'en\'eral \`a chaque \'el\'ement $x$ de U, on peut associer un ensemble
intuitif $\{y : y \in x\}$ . Cette application est m\^eme injective, par
l'axiome d'extentionalit\'e que U doit satisfaire. Mais on n'a pas de surjectivit\'e. A un ensemble intuitif de points de U ne correspond pas n\'ecessairement un point de U qui soit reli\'e par $\in^{\tt U}$ tr\`es exactement aux points de cet ensemble intuitif.
  
Un univers U est caract\'eris\'e par deux choses d'une part ses ordinaux, qui donnent sa hauteur, d'autre part l'op\'eration $\mathcal{P}$ de prise des parties, qui fournit l'\'epaisseur de U.

Un ordinal est un ensemble $a$ tel que si $x\in y\in a$, alors $x\in a$, et si $x\in a$ et $y\in a$, alors ou bien $x\in y$, ou bien $x=y$, ou bien $y\in x$. La relation $\in$ est donc un ordre total sur $a$, qui est un bon-ordre, c'est-\`a-dire est tel que tout sous-ensemble non vide a un plus petit \'el\'ement pour l'ordre. On obtient, \`a partir de l'ensemble vide (qui est un ordinal) tous les ordinaux par les op\'erations suivantes~:

Si $a$ est un ordinal, $a\cup\{a\}$ est un ordinal, appel\'e le successeur de $a$,
  
Si $X$ est un ensemble d'ordinaux, $a=\bigcup_{x\in X}\;x$ est un ordinal.

Par exemple, les premiers ordinaux sont $\emptyset$ (not\'e aussi $0$), puis son successeur $\{\emptyset\}$ , not\'e aussi 1, puis $\{\emptyset,\{\emptyset\}\}=2$, etc... Le premier ordinal autre que $\emptyset$ qui n'est pas successeur, et qui existe par l'axiome de l'infini, est l'ensemble, not\'e $\omega$ , des entiers $\omega = \{0,1,2,\ldots\}$.

Ces notions ont un sens dans U (qui est mod\`ele de la th\'eorie des ensembles).

\begin{center}
  \includegraphics[scale=1.4]{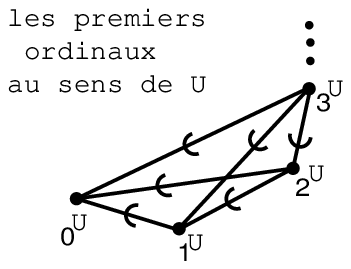}
\end{center}
\noindent A chaque entier correspond un entier $n^{\tt U}$ . Que dire de l'ensemble
intuitif $\{0^{\tt U},1^{\tt U},2^{\tt U},\ldots\}$~? On s'attendrait qu'il corresponde \`a l'ordinal $\omega^{\tt U}$. Pourtant il n'en est rien en g\'en\'eral cet ensemble intuitif peut tr\`es
bien ne correspondre \`a rien dans U, et l'ensemble $\omega^{\tt U}$ \^etre reli\'e \`a d'autres objets de U (qu'on appelle alors entiers non standards de U).

Une fois que les ordinaux de U sont connus, on d\'efinit pour chaque ordinal $\alpha$ de U un \'el\'ement ${\tt U}_\alpha$ de U par les clauses
$$
\begin{array}{l}
{\tt U}_0=\emptyset^{\tt U}\\
{\tt U}_\alpha=[\mathcal{P(\bigcup_{\beta\in\alpha}\;{\tt U}_\beta)]^{\tt U}}
\end{array}
$$
i.e. on obtient ${\tt U}_\alpha$ en prenant,au sens de U, la r\'eunion des ${\tt U}_\beta$ d\'ej\`a construits, puis l'ensemble des parties de cette r\'eunion.
Qu'une telle d\'efinition par r\'ecurrence soit possible est une propri\'et\'e fondamentale des ordinaux. Et on peut montrer qu'on obtient ainsi tout l'univers U, au sens suivant~: un \'el\'ement de U appara\^{\i}t comme \'el\'ement (au sens de U) de l'un des ${\tt U}_\alpha$.\\

\noindent{\large\bf4.	La construction de G\"{o}del}

Reprenons l'exemple de la th\'eorie des corps et du mod\`ele $\mathbb{R}$. Comment, \`a partir de ce mod\`ele, en construire un autre ? A priori, nous devons choisir une addition et une multiplication. Mais si on veut travailler \`a l'\'economie, on peut essayer de garder ces op\'erations, et ne changer que l'ensemble de base. Et notre souci d'\'economie nous am\`ene \`a construire le corps $\mathbb{Q}$. (En effet un sous-corps de $\mathbb{R}$ doit contenir 0 et 1 donc les entiers, donc les entiers relatifs et leurs inverses, donc tous les ration\-nels.)

Essayons la m\^eme id\'ee avec notre univers U. Nous n'allons pas
changer la relation d'appartenance $\in^{\tt U}$ (et donc mettre dans U', avec
un \'el\'ement $x$, tous ses \'el\'ements au sens de $\in^{\tt U}$). Les ordinaux de U \'etant
d\'efinis de fa\c{c}on tr\`es intrins\`eque \`a partir de $\in^{\tt U}$ , d\'ecidons aussi de conÂserver les ordinaux. Par contre, nous allons essayer de jouer sur l'op\'eration $\mathcal{P}$ pour diminuer au maximum l'"\'epaisseur" de l'univers.

Le probl\`eme est donc soit $x$ un \'el\'ement de U que nous avons d\'ecid\'e de conserver dans U' . Quelles parties de $x$ (au sens de U) devons-nous \'egalement conserver ? La r\'eponse est fournie par les axiomes de compr\'ehension~: nous devons mettre dans U' toutes les sous-parties de $x$ qui sont d\'efinies par des formules (avec param\`etres mis dans U').

Ceci donne l'id\'ee d'imiter la construction inductive des ${\tt U}_\alpha$ - qui fournit U - mais en d\'efinissant maintenant des $L_\alpha$ avec
$$
\begin{array}{ll}
L_{\alpha+1}=&{\{y\subset^{\tt U}L_\alpha\;|\;y=\{z\in^{\tt U}L_\alpha\;\&\;\varphi(z,a_1,\ldots,a_n)\}^{\tt U},} \\
&{\mbox{ pour une formule }\varphi\mbox{ et des param\`etres }a_1,\ldots,a_n\in^{\tt U}L_\alpha\}}.
\end{array}
$$

La d\'efinition qui pr\'ec\`ede est incorrecte~: l'ensemble d\'efini est intuitif, et doit \^etre remplac\'e par un \'el\'ement de U, ce qui est possible mais d\'elicat, une d\'efinition pr\'ecise des $L_\alpha$ peut \^etre trouv\'ee dans
le livre de J.L. Krivine [10]. Mais l'id\'ee de la construction de G\"{o}del
est bien celle que nous venons d'indiquer. Et il se trouve que l'ensemble
intuitif, not\'e g\'en\'eralement $L^{\tt U}$ , obtenu en conservant tous les \'el\'ements des $L_\alpha$ est bien un univers de la th\'eorie des ensembles, appel\'e univers
des ensembles constructibles de U. De plus, par sa construction, $L^{\tt U}$ est
le plus petit possible parmi les univers ayant m\^eme relation d'appartenance et m\^emes ordinaux. Et la construction des constructibles, refaite
dans $L^{\tt U}$ , fournit de nouveau	$L^{\tt U}$.

La construction de G\"{o}del permet donc de d\'emontrer la coh\'erence relative d'un \'enonc\'e de la Th\'eorie des Ensembles, l'\'enonc\'e V = L. Et par la m\^eme occasion, de toutes les cons\'equences de cet \'enonc\'e. Mais on peut prouver que parmi les cons\'equences de V = L figure l'hypoth\`ese du continu m\^eme g\'en\'eralis\'ee en "il n'y a pas de cardinalit\'e interm\'ediaire entre celle d'un ensemble A et celle de $\mathcal{P}(A)$", ainsi qu'un autre \'enonc\'e math\'ematique tr\`es discut\'e en son temps, l'axiome du choix.\\

\noindent{\large\bf5.	La construction de Cohen}

Pour construire un mod\`ele de la n\'egation de l'hypoth\`ese du continu, il faut d\'ej\`a savoir construire un mod\`ele de V $\neq$ L, et par ce qui pr\'ec\`ede, une diminution de l'univers est inefficace. La m\'ethode de Cohen consiste \`a l'augmenter, mais l\`a encore, \`a l'\'economie partant d'un univers U (que nous pouvons supposer satisfaire V = L), nous allons adjoindre un objet ext\'erieur $x$, pour fabriquer un nouvel univers U', mais de fa\c{c}on que la relation $\in^{\tt U' }$, restreinte aux \'el\'ements de U, soit la relation $\in^{\tt U}$
de U (un peu de la m\^eme fa\c{c}on que l'on passe de $\mathbb{R}$ \`a $\mathbb{C}$). D\'ecidons par exemple
de rajouter une nouvelle partie de $x$ de $\omega^{\tt U}$ . Nous allons certainement
devoir mettre dans U' d'autres objets, par exemple $\mathcal{P}(x)$. Comment s'asÂsurer que la plus petite collection U' ainsi obtenue est un mod\`ele de la th\'eorie des ensembles~? L'id\'ee de Cohen est de choisir $x$ de fa\c{c}on \`a perturber U le moins possible on va essayer de donner \`a $x$ un minimum de propri\'et\'es "particuli\`eres" si une propri\'et\'e des parties de $\omega^{\tt U}$ est vraie
dans U de "presque toutes" les parties de $\omega^{\tt U}$ , alors elle sera vraie dans
U' de $x$. La mise en forme math\'ematique de cette id\'ee est d\'elicate, mais fournit un mod\`ele U' dans lequel $x$ n'est pas constructible~: on prouve
ainsi la coh\'erence relative de V $\neq$ L. La construction de Cohen a l'avantage de pouvoir \^etre appliqu\'ee pour fournir un grand nombre de mod\`eles diff\'erents. Par exemple, en ajoutant beaucoup de nouvelles parties de $\omega^{\tt U}$
\`a U, on peut obtenir un mod\`ele U' dans lequel $\mathcal{P}(\omega^{\tt U})^{\tt U}$ a une
cardinalit\'e interm\'ediaire entre celle de $\omega^{\tt U}$ et celle de $\mathcal{P}(\omega^{\tt U})^{\tt U'}$ , c'est-\`a- dire un mod\`ele dans lequel la n\'egation de l'hypoth\`ese du continu est satisfaite. L'\'etude des mod\`eles qui sont produits par la m\'ethode de Cohen a permis, durant ces vingt derni\`eres ann\'ees, de montrer qu'un grand nombre d'\'enonc\'es math\'ematiques sont ind\'ecidables, touchant (de mani\`ere tr\`es in\'egale) de nombreuses branches des math\'ematiques.\\

\noindent{\large\bf6.	Et l'hypoth\`ese du continu~?}

La r\'eponse apport\'ee par G\"{o}del et Cohen au probl\`eme de Cantor est d\'efinitive. Mais elle n'est pas n\'ecessairement convaincante. G\"{o}del lui-m\^eme n'en \'etait pas du tout satisfait cette r\'eponse indique en effet plut\^ot les limites de la formalisation, et une r\'eponse par vrai ou faux est envisageable, \`a condition d'augmenter les axiomes de la th\'eorie des ensembles. On pourrait bien s\^ur ajouter comme axiome l'hypoth\`ese du continu elle-m\^eme, ou sa n\'egation. Mais que choisir~? Une autre possibilit\'e serait de prendre comme nouvel axiome l'\'enonc\'e de G\"{o}del  V = L. C'est un \'enonc\'e qui a de nombreuses cons\'equences int\'eressantes, mais il a le d\'efaut de "diminuer" l'univers d'etude, ce qui, ne serait-ce que d'un point de vue m\'ethodologique, n'est pas du tout satisfaisant. On peut imaginer que l'on propose un jour un nouvel axiome, \`a la fois naturel et intuitivement vrai, qui r\'esolve l'hypoth\`ese du continu, et que la pratique math\'ematique d'alors en fasse l'une des "r\'ealit\'es math\'ematiques" indiscut\'ees. En attendant, on se trouve place devant une pluralit\'e de choix (et ce, pas seulement pour l'hypoth\`ese du continu), avec laquelle nous devons bien vivre.

\newpage
\noindent{\LARGE\bf Sur la calculabilit\'e effective, exemples}

\noindent Nicolas Bouleau\\

\noindent{\large\bf I.	Syst\`emes formels}

Le syst\`eme formel le plus simple est le "syst\`eme des allumettes".
il	est constitu\'e d'un seul signe I et ses assemblages sont des groupes finis de ce signe r\'ep\'et\'e I, II, III, IIII, IIIII,...

L'arithm\'etique des entiers naturels d\'ecrit des lois de ce syst\`eme formel, elle-m\^eme formalis\'ee		suivant un syst\`eme
formel~: l'arithm\'etique de Peano.

C'est une th\'eorie du 1er ordre (cf. 1er expos\'e).

\noindent Langage du 1er ordre

	a) constante 0
	
	b) fonction unaire S (successeur),
fonctions binaires + et .
	
c)	pr\'edicat binaire <
	
\noindent Axiomes.	
	
	P1.   $Sx\neq 0$
	
	P2.	$Sx=	Sy\Rightarrow x=y$
	
	P3.	$x+0	=x$
	
	P4.	$x+Sy=S(x+y)$
	
	P5.	$x.0=0$
	
	P6.	$x.Sy=(x.y) + x$
	
	P7.   $\neg(x<0)$
	
	P8.	$x <Sy\Leftrightarrow x<y\vee x=y$

\noindent et le sch\'ema d'axiomes suivant

P9. Si $A(x)$ est une formule admettant $x$ pour variable libre 
	
	$A(0) \;\&\; \forall y(A(y)\Rightarrow A(S(y)) \Rightarrow \forall x A(x)$

Est-ce que les th\'eor\`emes de cette th\'eorie sont vrais~? Si l'on regarde ce que signifie tel th\'eoÂr\`eme pour notre syst\`eme des allumettes et si on effectue sur ces derni\`eres les op\'erations indiqu\'ees on constate que \c{c}a marche. Par ses cons\'equences dans toutes les math\'ematiques, dans la physique et en informatique on peut dire que cette th\'eorie est de loin celle qui a \'et\'e le plus s\'ev\`erement test\'ee.
	
Une fois convaincu de cela il reste cependant deux questions importantes

1$^0$)	A-t-on en recueillant tous les th\'eor\`emes de cette th\'eorie, toutes les propri\'et\'es vraies des allumettes qui peuvent s'exprimer dans le langage utilis\'e ?

2$^0$)	Existe-t-il un moyen m\'ecanique de savoir si une propri\'et\'e des allumettes exprim\'ee en ce langage est cons\'equence des axiomes ou non~?

Je vais tenter de faire comprendre intuitivement comment dans les ann\'ees trente, les logiciens (Turing, G\"{o}del, Church, Kleene, etc.) ont pu r\'epondre \`a ces questions.

Prenons un syst\`eme formel tel que celui de Peano ou celui de Zermelo-Fraenkel. Bien que nous ayons plusieurs r\`egles de d\'eduction \`a notre disposition, il est possible pour \'ecrire nos th\'eor\`emes, de proc\'eder de fa\c{c}on syst\'ematique.

En effet, on peut num\'eroter les symboles de notre langage puis num\'eroter les termes et les formules et enfin puisque les d\'emonstrations sont des suites finies de formules, chacune reli\'ee aux pr\'ec\'edentes par un nombre fini de r\`egles, on peut num\'eroter les d\'emonstrations\footnote{Voir les expos\'es de J.-Y. Girard.}.

Alors si nous prenons les d\'emonstrations par ordre de num\'eros croissants, nous passons ainsi en revue tous les th\'eor\`emes de notre syst\`eme.

La chose importante \`a comprendre est que si nous proc\'edons de cette mani\`ere syst\'ematique, les th\'eor\`emes que nous obtenons successivement sont des \'enonc\'es de longueur tr\`es variable. Nous n'obtenons pas les th\'eor\`emes par ordre de num\'eros croissants, il se peut en effet qu'une d\'emonstration tr\`es longue m\`ene \`a un th\'eor\`eme assez court pour lequel il n'y ait pas de d\'emonstration plus courte. De telle sorte que si, \`a l'inverse, nous nous donnons un \'enonc\'e a priori et si nous cherchons \`a savoir si c'est un th\'eor\`eme de notre th\'eorie, { il n'y a pas de borne sup\'erieure \`a la longueur
des d\'emons\-trations \`a investiguer pour savoir si notre \'enonc\'e est d\'eductible ou non de nos axiomes}.

Les logiciens ont alors d\'emontr\'e les r\'esultats suivants

Si le syst\`eme formel est suffisamment riche (grosso-modo s'il est au moins aussi riche que l'arithm\'etique de Peano), et s'il n'est pas contradictoire (c'est-\`a-dire s'il existe des \'enonc\'es qui ne sont pas cons\'equences des axiomes) alors

1) Il n'existe aucun algorithme (c'est-\`a-dire de m\'ethode programmable sur un ordinateur m\^eme de m\'emoire de capacit\'e infinie) permettant de dire en un temps fini si un \'enonc\'e est un th\'eor\`eme ou non.

2) Il existe m\^eme des \'enonc\'es tels que ni eux ni leur n\'egation ne sont des th\'eor\`emes.

La premi\`ere propri\'et\'e s'exprime en logique par le fait que l'ensemble des th\'eor\`emes est un ensemble recursivement \'enum\'erable non r\'ecursif, ou encore en disant que pour un \'enonc\'e, le fait d'\^etre un th\'eor\`eme est une propri\'et\'e semi-d\'efinie positive.

Pour faire comprendre cela intuitivement prenons une analogie et consi\-d\'erons un ouvrier charg\'e d'examiner un lot infini de pommes passant devant lui sur un tapis roulant. Si d'aventure il rencontre une pomme v\'ereuse alors il est s\^ur qu'il y a des pommes v\'ereuses dans le lot, mais si nous imaginons que l'ouvrier ne rencontre jamais de pommes v\'ereuses alors l'ouvrier ne saura jamais s'il y a des pommes v\'ereuses ou non dans le lot. C'est moins trivial qu'il n'y para\^{\i}t si l'on pose maintenant la question existe-t-il un groupe de 10 z\'eros \`a la suite dans les d\'ecimales de $\sqrt{2}$~? Si l'on d\'ecouvre une telle suite, la r\'eponse est positive, si on n'en trouve pas on ne sait pas r\'epondre \`a la question.

A moins que l'on ne dispose par bonheur d'une theorie dont un des th\'eo\-r\`emes dit justement "Il existe dans les d\'ecimales de $\sqrt{2}$ un groupe de 10 z\'eros \`a la suite". On ne conna\^{\i}t pas actuellement de tel th\'eoÂr\`eme et on est donc oblig\'e de passer en revue les th\'eor\`emes et si l'on ne trouve jamais de tel th\'eor\`eme on est dans la m\^eme situation qu'avant.

Ce type de r\'eflexion se rattache historiquement \`a l'intuitionisme que je ne puis aborder maintenant.
	
La deuxi\`eme propri\'et\'e s'\'enonce en disant que le syst\`eme est incomplet. Historiquement le premier \'enonc\'e d'arithm\'etique ni prouvable ni r\'efutable fut d\'ecouvert par Rosser am\'eliorant un travail pr\'ealable de G\"{o}del. A noter que quoique formellement ind\'ecidable cet \'enonc\'e est intuitivement vrai. Il en est de m\^eme pour la traduction dans le langage de l'arithm\'etique de l'assertion "L'arithm\'etique est non contradictoire". Cet \'enonc\'e n'\'etant pas r\'efutable on peut \'evidemment l'ajouter comme nouvel axiome sans risque nouveau de contradiction. Cependant la nouvelle th\'eorie obtenue aura les m\^emes tares d'ind\'ecidabilit\'e et d'incompl\'etude que l'ancienne.\\

\noindent{\large\bf II.	Fonctions r\'ecursives}

On peut d\'emontrer dans ZF le th\'eor\`eme suivant\\

\noindent{\bf Th\'eor\`eme.} {\it Pour presque tout r\'eel $\theta>1$ pour la mesure de Lebesgue, la suite $a_n=\mbox{ reste de } \theta^n\mbox{ modulo }1$, est \'equir\'epartie sur $[0,1]$,
	
\noindent C'est-\`a-dire si $[a,b]\subset[0,1]$

$$\lim_{N\uparrow\infty}\frac{\mbox{nombre des }a_n\in[a,b]\mbox{ pour }n\leq N}{N}=b-a.
	$$}
	
Ce qui est chagrinant par ailleurs c'est qu'on ne sait citer concr\`etement aucun nombre $\theta>1$ qui v\'erifie la propri\'et\'e.

On a un th\'eor\`eme d'existence mais ce th\'eor\`eme ne permet pas de construire effectivement un objet dont on affirme l'existence.

Un nombre $\theta$ r\'eel en base 2 est une application de $\mathbb{N}$ dans $\{0,1\}$, quand pourra-t-on dire qu'une application de $\mathbb{N}$ dans lui-m\^eme est effectivement calculable~? Une r\'eponse \`a cette question est donn\'ee par la th\'eorie des fonctions r\'ecursives.

Comme la notion de fonction "effectivement calculable" est intuitive et un peu floue, les math\'ematiciens lui pr\'ef\`erent une notion plus pr\'ecise (a priori plus restreinte) de fonction recursive.

Les fonctions r\'ecursives d'entiers \`a valeurs enti\`eres sont d\'efinies de la fa\c{c}on suivante

1) les projections de $\mathbb{N}^n$ dans $\mathbb{N}$ sont r\'ecursives.
	$$(a_1,\ldots,a_i,\ldots,a_n)\mapsto a_i$$
	
	2)	les fonctions de 2 variables
	$$\begin{array}{l}
	.+.\\
	.\times.\\
	1_{\{.<.\}}
	\end{array}
	$$
	 sont r\'ecursives\footnote{La notation $1_A$ d\'esigne la fonction indicatrice de l'ensemble $A$~: elle vaut 1 au point $x$ si $x$ est dans $A$, z\'ero sinon.}

3) si $G, H_1,\ldots,H_n$ sont r\'ecursives la fonction compos\'ee
	$G(H_1,\ldots,H_n )$ est r\'ecursive.

4) si $G(a_0,a_1,\ldots,a_n)$ est r\'ecursive et si $G$ v\'erifie 
	$$\forall a_1,\ldots, \forall a_n,\exists a : G(a,a_1,\ldots,a_n)=0$$
		la fonction
	$$F(a_1,\ldots,a_n) = \mbox{ le plus petit } a \mbox{ tel que }
	 G(a,a_1,\ldots,	a_n) = 0$$
est recursive.

Je vous renvoie aux livres sur la calculabilit\'e effective o\`u sont d\'evelopp\'es de proche en proche des exemples de fonctions r\'ecursives pour vous convaincre de ce qu'on appelle la \underline{th\`ese de Church}, c'est-\`a-dire du fait que les fonctions effectivement calculables sont exactement les fonctions r\'ecursives.

Nous dirons qu'un sous-ensemble de $\mathbb{N}$ est r\'ecursif si sa fonction indicatrice est r\'ecursive.

Il est clair que par application successives des r\`egles 1) 2) 3) 4) on obtient un ensemble d\'enombrable de fonctions~: l'ensemble des fonctions r\'ecursives est d\'enombrable. Au demeurant la num\'erotation des foncÂtions r\'ecursives ne peut elle-m\^eme \^etre r\'ecursive, par un argument diagonal simple~: s'il existe une fonction r\'ecursive $h(m,n)$ telle que lorsque $m$ varie les fonctions $n \mapsto h(m,n)$ soient toutes les fonctions r\'ecursives d'une variable, la fonction $n\mapsto h(n,n) + 1$ serait r\'ecursive et ne serait pas dans la liste.

La particularit\'e d'un ensemble non r\'ecursif $A$ est qu'il n'existe aucun algorithme permettant de r\'epondre \`a la question
$$n\in A\mbox{ ?}$$
c'est-\`a-dire de m\'ethode de calcul qui appliqu\'ee \`a $n$ r\'epondra en un nombre fini d'\'etape \`a la question
$n\in A$~?

Les fonctions r\'ecursives permettent de pr\'eciser la notion de "propri\'et\'e semi-d\'efinie positive" gr\^ace \`a la notion rigoureuse d'"ensemble r\'ecursivement \'enum\'erable"~: si $n \mapsto f(n)$ est r\'ecursive, l'image de la fonction $f$, c'est\-\`a-dire l'ensemble $\{n\;|\;\exists m \; f(m)=n\}$ est un ensemble r\'ecursivement \'enum\'erable. On peut montrer que c'est une classe d'ensembles \underline{plus vaste} que celle des ensembles r\'ecursifs. Si $A$ est un ensemble r\'ecursivement \'enum\'erable non r\'ecursif, alors la question $n\in A$ ? est semi-d\'efinie positive~: si $n\in A$ un calcul effectif nous le dira en un temps fini, mais si $n\notin A$, il faut passer en revue tous les entiers pour s'assurer qu'aucun ne donne $n$ comme image par $f$.

Ainsi l'application d'un quantificateur existentiel devant une propri\'et\'e r\'ecursive $P(m,n)$ donnant $\exists m P(m,n)$ conduit \`a des ensembles de plus grande complexit\'e.

Le ph\'enom\`ene est analogue au fait que la projection sur $\mathbb{R}$ d'un
bor\'elien de $\mathbb{R}^2$ n'est pas en g\'en\'eral un bor\'elien (mais un ensemble analytique comme l'a montr\'e Souslin, contrairement \`a ce que Lebesgue avait cru pouvoir affirmer).

En revanche, on peut montrer que si $A$ et son compl\'ementaire sont r\'ecur\-sivement \'enum\'erables alors ils sont en fait r\'ecursifs.

Un des r\'esultats les plus c\'el\`ebres et li\'e \`a ce que j'\'evoquais au d\'ebut de l'expos\'e est le suivant\\

\noindent{\bf Th\'eor\`eme de Church.}

{\sl Consid\'erons une num\'erotation effective des th\'eor\`emes d'une th\'eorie non contradictoire qui contienne l'arithm\'etique alors l'ensemble des th\'eor\`emes est non r\'ecursif.}\\

\noindent{\large\bf III. Exemples}

Nous allons maintenant donner quelques exemples de probl\`emes r\'ecur\-sivement insolubles.\\

\noindent a) Dixi\`eme probl\`eme de Hilbert.

Hilbert au congr\`es international des math\'ematiciens \`a Paris en 1900 proposait \`a la post\'erit\'e une liste de 23 probl\`emes dont la solution \`a son avis risquait de faire avancer les math\'ematiques.

Le 10e probl\`eme a trait aux \'equations diophantiennes, c'est-\`a-dire la r\'esolution en nombres entiers d'\'equations polynomiales \`a coefficients entiers par exemple du type
$$4x^2-xy^3+7z^5+1=0.$$
 La question pos\'ee par Hilbert \'etait la suivante

{\it Existe-t-il un algorithme permettant de dire si une \'equation diophantienne admet une solution ou non ?}

 Hilbert s'int\'eressait \`a ce probl\`eme et attachait de l'importance \`a la notion d'algorithme d\`es cette \'epoque.

Il s'agit donc de trouver une m\'ethode permettant pour tout polyn\^ome \`a coefficients entiers, $P(x_1,\ldots,x_n)$
de dire s'il existe $n$ nombres entiers $a_1,\ldots,     a_n$ , tels que
$P(a_1,\ldots,a_n)=0$  ou non.
	 
La r\'eponse a \'et\'e r\'esolue par la n\'egative par Yu. V. Matiyassevic en 1970, prolongeant les travaux de Putnam, Davis et Julia Robinson.

La d\'emonstration permet de montrer le r\'esultat plus pr\'ecis suivant~:\\

\noindent {\it Il existe un polyn\^ome de 10 variables \`a coefficients entiers
$$U(a,x_1,\ldots,x_9)$$
tel qu'il n'existe aucun algorithme pour dire si les \'equations
	$$n\mbox{ fix\'e }\quad U(n,x_1,\ldots,x_9)=0$$
ont des solutions enti\`eres ou non.}\\

Et comme sous-produit de la d\'emonstration on obtient aussi l'existence d'un polyn\^ome de 10 variables \`a coefficients entiers dont l'ensemble des valeurs prises lorsque les variables parcourent les entiers est exactement l'ensemble des nombres premiers.
	 
Enfin, en traduisant le fait qu'on peut d\'emontrer dans ZF la non-contra\-diction du syst\`eme de Peano, on obtient le r\'esultat suivant~:\\

\noindent{\it Il	existe un polyn\^ome de plusieurs variables \`a coefficients
entiers $P$ (que l'on pourrait \'ecrire explicitement) pour lequel c'est un th\'eor\`eme de ZF que l'\'equation $P=0$
n'a aucune solution enti\`ere, mais pour lequel cette assertion n'est pas d\'emontrable dans l'arithm\'etique de Peano.}\\

La th\'eorie ZF am\'eliore donc notre connaissance de l'arithm\'etique, nous reviendrons sur ce point en conclusion.\\

\noindent b)	Le probl\`eme des mots pour un semi-groupe.

On consid\`ere le semi-groupe libre \`a 2 g\'en\'erateurs c'est-\`a-dire les mots form\'es de 2 lettres $a$ et $b$
$$a, b, aa, ab, aba,\ldots, abbabab,\ldots$$
muni de la loi de composition qui est la juxtaposition sans parenth\`eses donc associative.

Nous d\'esignons les mots par des lettres majuscules.
On d\'efinit alors une relation d'\'equivalence sur l'ensemble des mots de la fa\c{c}on suivante~:

On prend deux mots $A = abbabab$ et $B = aaba$ par exemple et on pose
$$A\equiv B \quad\mbox{ (c'est l'axiome)}$$
et on impose les r\`egles suivantes entre mots
	 $$
	 \begin{array}{rl}
	 C\equiv C'&\Rightarrow CD\equiv C'D\\
	 C\equiv C'&\Rightarrow EC\equiv EC'\\
	 FGH\equiv FG'H&\Rightarrow G\equiv G'
	 \end{array}
	 $$
Alors on peut montrer que:\\

\noindent{\it Il	existe des axiomes en nombre finis
	$$A_1\equiv B_1,\ldots,A_n\equiv B_n$$
(qu'on peut \'ecrire effectivement) tels qu'il n'existe aucun algorithme permettant de r\'esoudre la question
$$E\equiv F\mbox{ ?}$$
entre deux mots quelconques.}\\
	 
\noindent c)	Linguistique formelle.

Je renvoie pour les d\'etails \`a
Davis M.,	{"Unsolvable problems"},in {\it Handbook of Mathematical Logic} Barwise editor North Holland 1978.\\

\noindent{\it Il	n'y a pas d'algorithmes pour dire si deux grammaires sans contexte permettent une phrase commune ou non.}\\

Une grammaire est dite ambigu\"{e} si elle permet des phrases ambigu\"{e}s c'est-\`a-dire dont la structure grammaticale peut \^etre de deux types diff\'erents.
Par exemple la phrase

"Pourquoi crois-tu que Walesa est en prison ?"

\noindent est ambigu\"{e}.
Si pourquoi se rapporte au groupe "Walesa est en prison" 
elle recevra une r\'eponse du type~:

\noindent R\'eponse : parce que les autorit\'es ont jug\'e qu'il \'etait dangereux. Si pourquoi se rapporte au groupe "crois-tu"
elle recevra une r\'eponse du type
	 
\noindent R\'eponse : parce que c'est ce que dit la presse.

\noindent{\it Il n'y a pas d'algorithme permettant de dire si une grammaire formelle est ambigu\"{e} ou non.}\\

\noindent d)	Autre exemple.

Il n'y a pas d'algorithme permettant de dire si deux vari\'et\'es de dimensions 4 sont hom\'eomorphes ou non.\\

\noindent e)	En revanche des th\'eories plus faibles que l'arithm\'etique (dans lesquelles on ne peut pas faire la th\'eorie des fonctions r\'ecursives ou celle des machines de Turing) peuvent \^etre d\'ecidables, c'est le cas pour la th\'eoÂrie de l'addition de Presburger. Cette th\'eorie a le m\^eme langage que l'arithm\'etique de Peano except\'e la fonction binaire $.\times.$, et les m\^emes axiomes except\'es P5 et P6.

Dans cette th\'eorie l'ensemble des th\'eor\`emes est r\'ecursif. Il existe un algorithme permettant de dire en un temps fini si une formule
du langage est un th\'eor\`eme ou non. Un \'enonc\'e est soit prouvable soit r\'efutable (i.e. sa n\'egation est prouvable) : la th\'eorie est compl\`ete. Il existe une m\'ethode finitiste permettant de montrer que cette th\'eorie est non contradictoire.\\

\noindent{\large\bf IV.	Quelques commentaires philosophiques : l'exemple et l'argument}

Il est toujours assez hasardeux de tirer des cons\'equences philosophiques de travaux scientifiques, d'une part ces commentaires sont parfois de peu d'utilit\'e pour la recherche scientifique elle-m\^eme, d'autre part c'est leur nature, semble-t-il, de vieillir assez rapidement.

Si je m'y risque ici c'est pour critiquer quelque peu des conceptions philosophiques pr\'ecis\'ement fond\'ees sur une certaine image de la science et des math\'ematiques en particulier.

Une question qui a fait couler beaucoup d'encre depuis le 17e si\`ecle est de concilier la f\'econdit\'e du raisonnement math\'ematique avec sa nature rigoureuse. Les r\'ef\'erences sont abondantes, on peut citer, par exemple, Pascal, Kant, Poincar\'e, les positivistes logiques (Carnap, etc.).

La d\'ecouverte de l'ind\'ecidabilit\'e \'eclaire ce probl\`eme d'un jour nouveau.
Il faut dire \`a ce sujet que si, de 1930 \`a la guerre, les positivistes logiques avaient bien "assimil\'e" la formalisation des math\'ematiques d\'ecouverte 20 ans plus t\^ot, il n'en va pas de m\^eme pour les questions relatives \`a l'ind\'eci\-dabilit\'e qui fut \'elucid\'ee pour une large part par des math\'ematiciens (G\"{o}del, Tarski) li\'es au cercle de Vienne lui-m\^eme. Des le\c{c}ons de ces travaux ne seront tir\'es qu'apr\`es la guerre chez Putnam, Quine, voire Chomsky, etc.\\

\noindent A.	Lorsqu'un m\'ecanisme combinatoire est assez simple, on peut dire que tous les assemblages qu'il engendre sont "contenus" dans ses hypoth\`eses et ses r\`egles de formation, il ne r\'eserve en v\'erit\'e aucune surprise et une connaissance rationnelle peut l'embrasser globalement.

Au contraire lorsque le m\'ecanisme atteint un certain niveau de complexit\'e (caract\'eris\'e par le fait que ses "outputs" forment un ensemble
non r\'ecursif) alors il n'existe aucune m\'ethode de connaissance certaine des r\'esultats de ces mecanismes autre que l'exp\'erimentation du m\'ecanisme lui-m\^eme\footnote{ou d'un autre mecanisme, li\'e au pr\'ec\'edent ou le contenant, Cette nuance est \'epist\'emologiquement importante dans la mesure o\`u les syst\`emes formels que sont AP ou ZF sont consid\'er\'es par tous les math\'ematiciens comme absolument fiables. Cette conviction solide ne s'appuie \'evidemment pas sur une preuve formelle, sur quoi se fonde-t-elle ? Il serait long de l'analyser. En tout cas ne pas croire \`a la coh\'erence de ZF risque fort d'\^etre une attitude math\'ematiquement st\'erile. Le bon choix est de faire le pari, c'est l'argument pascalien avec cette diff\'erence qu'ici on y gagne m\^eme en ce bas monde.}.

On ne sait pas quels sont les th\'eor\`emes de l'arithm\'etique, on ne conna\^{\i}t que des exemples de tels th\'eor\`emes.

On peut dire en quelque sorte que le syst\`eme quoique purement m\'ecanique est d'une complexit\'e suffisante pour que son d\'eveloppement engendre en permanence des surprises.

Dans la dualit\'e classique entre d\'eterminisme et hasard un troisi\`eme terme intervient qui est l'ind\'ecidable. Et pour reprendre l'image tant comment\'ee de Laplace : on pourrait dire que, quand bien m\^eme seraient connues toutes les positions des mol\'ecules de l'univers et les champs de forces agissant sur elles que, sans qu'aucun hasard n'intervienne, leurs combinaisons ne cesseraient de nous r\'eserver des surprises.

Plusieurs auteurs contemporains dont Edgar Morin par exemple ont vu l'importance de la complexit\'e dans les probl\`emes \'epist\'emologiques (de la biologie surtout~: ontogen\`ese et phylogen\`ese). H\'elas, cette litt\'erature est assez confuse et les liens avec les syst\`emes formels ne sont pas clairement d\'egag\'es de tout le fourbi thermo-dynamico-informationnel. Cependant d'importants travaux r\'ecents (Kolmogorov, Rabin, Blum, Constable, Chaitin, etc) tendent \`a pr\'eciser cette notion de complexit\'e.\\

\noindent B.	Revenons aux math\'ematiques, ce que nous avons dit \`a l'instant pourrait faire conclure qu'\'epist\'emologiquement c'est l'exemple qui est le
plus fort. Cependant les choses ne sont, l\`a encore, pas si simples. Nous
ne savons pas manier les syst\`emes formels en eux-m\^emes sans attribuer
une signification aux symboles\footnote{Il n'est pour s'en convaincre que d'\'ecrire et d'essayer de d\'emontrer sans aucune abr\'eviation dans le langage de ZF que la limite uniforme d'une suite de fonctions r\'eelles continues est continue.}. Nous r\'eussissons mieux \`a faire de la
th\'eorie que de la combinatoire. C'est l\`a que l'argument reprend sa place.

Plut\^ot que de chercher vainement une d\'emonstration de telle conjecture arithm\'etique, on peut avec peut-\^etre plus de succ\`es tenter de construire une th\'eorie math\'ematique plus puissante que l'arithm\'etique dans laquelle on a une bonne confiance et de regarder les propositions arithm\'etiques qu'elle permet de d\'emontrer (cf. III in fine). C'est une des raisons qui poussent certains logiciens contemporains \`a rechercher des axiomes nouveaux \`a adjoindre \`a ZF, souvent ces axiomes sont d'un m\^eme type : axiomes de l'infini. Ils postulent l'existence d'un cardinal tr\`es grand. Ces axiomes permettent de d\'emontrer des th\'eor\`emes d'arithm\'etique qu'on ne sait pas d\'emontrer dans le syst\`eme de Peano ni m\^eme dans
ZF.

Dans cette course aux cardinaux tr\`es grands le risque \'evident est de viser trop haut et d'aboutir \`a une th\'eorie contradictoire.

On peut dire en quelque sorte que les th\'eories les plus fructueuses  semblent \^etre parmi celles qui prennent le plus de risques vis-\`a-vis du probl\`eme de non contradiction.
On trouve une situation analogue \`a celle que d\'ecrit Popper \`a propos des sciences de la nature lorsqu'il dit qu'entre deux th\'eories la plus int\'eressante est celle qui prend le plus de risques vis-a-vis des sanctions de l'exp\'erience.

Mais la discussion ne s'arr\^ete pas l\`a. Les th\'eories ne se laissent pas classer en ordre lin\'eaire. Il y a des fa\c{c}ons non comparables d'exprimer qu'un cardinal est grand. Et de m\^eme, pour les autres sciences \underline{l'initiative} revient au chercheur et c'est \`a lui d'assumer les risques engendr\'es par son activit\'e.
	 
\newpage

\noindent{\LARGE\bf Bibliographie}\\
\vspace{.5cm}

\noindent{\bf I -}	Un livre de vulgarisation sur les math\'ematiques et leurs probl\`emes historiques, philosophiques et \'epist\'emologiques

\noindent[1] R. HERSH \& Ph. J. DAVIS
{\it The Mathematical Experience},
Birkh\"{a}user Boston 1980.\\

\noindent {\bf Il	- Livres introductifs}

\noindent[2] M. COMBES	{\it Fondements des math\'ematiques},
P.U.F.	coll. sup. 1971.

\noindent[3] B.C. LYNDON {\it Notes on logic}
Von Nostrand 1964.\\

\noindent {\bf III	- Livres g\'en\'eraux}

\noindent[4] S.C. KLEENE{ \it Introduction to Metamathematics},
North Holland 1952.
(Le livre du m\^eme auteur {\it Mathematical logic}, Wiley 1967 - traduit en fran\c{c}ais chez A. Colin est moins riche).

\noindent[5] SHOENFIELD {\it Mathematical Logic}
Addison Wesley 1967.

\noindent signalons \'egalement l'ouvrage collectif

\noindent[6] BARWISE \& al., {\it Handbook of Mathematical Logic}
North Holland 1978.

\noindent o\`u l'on trouvera des bibliographies d\'etaill\'ees.\\

\noindent {\bf IV	Grandes branches de la logique}

\noindent A.	Th\'eorie de la r\'ecursivit\'e

Les r\'ef\'erences sont tr\`es nombreuses, un classique est le livre de

\noindent[7] M. DAVIS,	{\it Computability and unsolvability}
Mc Graw Hill 1958.

voir \'egalement [5] et les bibliographies de [6].\\

\noindent B.	Th\'eorie des mod\`eles

\noindent[8] CHANG \& KEISLER, {\it Model Theory}
North Holland 1973.

\noindent[9] BELL \& SLOMSON, {\it Models and Ultraproducts,  an introduction} North Holland 1974.\\

\noindent C.	Th\'eorie des ensembles

\noindent[10] J.L. KRIVINE, {\it Th\'eorie axiomatique des Ensembles}
P.U.F.	1972.\\

\noindent D.	Th\'eorie de la demonstration

voir [6] et [12] [13] [14] ci-dessous.\\

\noindent {\bf V	- R\'ef\'erences relatives aux cinq conf\'erences}

\noindent A.	L'expos\'e sur la formalisation des math\'ematiques utilise abondamment l'\'etude historique de

\noindent[11] M. GUILLAUME "Axiomatique et logique"
in {\it Abr\'eg\'e d'histoire des mathema-tiques} Tome II, J. Dieudonn\'e Hermann 1978.\\

\noindent B.	A propos des expos\'es de J.Y. GIRARD on pourra consulter

\noindent[12] K. G\"{O}DEL	"Ueber formale unentscheidbare S\"{a}tze der Principia Mathematica und verwandter Systeme" Monatshefte
f\"{u}r Math. und Physik. 38, pp 173 - 198,	1931.

\noindent[13] C. GENTZEN	"Die Widerspruchfreiheits beweise der reine Zahltheorie"
{\it Math. Ann.} 112 p 493 - 595 (1936)

\noindent[14] G. KREISEL "A survey of Proof theory"
		{\it J.	of symb. logic}, vol. 33 (1968).

"A survey of Proof theory II"
{\it Proc. of the 2d scandinavian logic
symp.} Feustad ed. North Holland (1971).\\

\noindent C.	Au sujet de l'expos\'e d'Alain LOUVEAU sur l'hypoth\`ese du continu on pourra consulter les textes originaux suivants

\noindent[15] K. G\"{O}DEL,	"The consistency of the axiom of choice and the generalized continuum hypothesis".
{\it Proc. Nat. Acad. Sci. USA}, 24, pp 556-557
(1938).
 
\noindent[16] P. COHEN	{\it Set Theory and the Continuum Hypothesis} Benjamin 1966.

\noindent ainsi que les cours sur la th\'eorie du forcing~:

\noindent[17] J.L. KRIVINE	{\it Th\'eorie des ensembles}
Cours du 3e cycle multigraphie 1970.

\noindent[18] S. GRIGORIEFF \& J. STERN
{\it Th\'eorie et pratique du forcing}
Cours du 3e cycle multigraphi\'e
Universit\'e Paris VII.\\

\noindent D.	Sur la calculabilit\'e effective on pourra consulter

\noindent DAVIS, MATIYASSEVITCH, J. ROBINSON
{\it Procedings of symp. in pure math.} vol. 28 pp 323-378 (1976).

\noindent M.	MARGENSTERN	"Le th\'eor\`eme de Matiyassevitch et r\'esultats connexes"
			                 preprint.

\noindent un aper\c{c}u des recherches r\'ecentes sur la complexit\'e est donn\'e dans

\noindent J.E. HOPCROFT, J.D. ULLMAN
{\it Introduction to Automata Theory, Languages and
		Computation}
Addison - Wesley (1979).

}

\end{document}